\documentclass[letterpaper, 10 pt, conference]{IEEEtran}
\usepackage{cite}
\usepackage{afterpage}
\usepackage{algorithm}
\usepackage[noend]{algpseudocode}
\usepackage{pifont}
\usepackage{wrapfig}
\usepackage{amsmath,amssymb}
\usepackage{bm}
\usepackage[authormarkup=none,final]{changes}
\usepackage{epigraph} 
\usepackage{float}
\usepackage{graphicx}
\usepackage{pdflscape}
\usepackage{lettrine}
\usepackage{mathtools}
\usepackage{mathrsfs}
\usepackage{multirow}
\usepackage{hyperref}
\usepackage[numbers]{natbib}
\usepackage[percent]{overpic}
\usepackage[font={footnotesize,it}]{caption}
\usepackage{subcaption}
\usepackage{tikz}
\usepackage{upgreek}
\usepackage{xcolor}
\usepackage{xfrac}
\usepackage{siunitx}
\usepackage{cuted}
\usepackage{tabularx}
\usepackage[cal=dutchcal]{mathalfa} 

\def\BState{\State\hskip-\ALG@thistlm}
\newcommand{\MATLAB}{\textsc{Matlab}\textregistered\:}

\newcommand{\matrixstyle}[1]{\mathrm{#1}}
\newcommand{\vectorstyle}[1]{\boldsymbol{\mathbf{#1}}}
\renewcommand{\baselinestretch}{0.95}
\DeclareMathOperator*{\Par}{Par}

\def\IEEEtitletopspace{18pt}
\usepackage[top=54pt,bottom=54pt,left=54pt,right=54pt]{geometry}
\IEEEoverridecommandlockouts

\begin{document}
	
\title{Introducing DAIMYO: a first-time-right dynamic design architecture and its application to tail-sitter UAS development}

\author{Jolan~Wauters,~Tom~Lefebvre,~Joris~Degroote,~Ivo~Couckuyt~and~Guillaume~Crevecoeur
	\thanks{J. Wauters, T. Lefebvre and G. Crevecoeur are with the Dynamic Design Lab (D2Lab) of the Department of Electromechanical, Systems and Metal Engineering (ESME), Ghent University, B-9052 Ghent, Belgium e-mail: \{jolan.wauters, tom.lefebvre, guillaume.crevecoeur\}@ugent.be.}%
	\thanks{J. Degroote is with the Fluid Mechanics group of the Department of Electromechanical, Systems and Metal Engineering, Ghent University, B-9052 Ghent, Belgium e-mail: joris.degroote@ugent.be.}%
	\thanks{I. Couckuyt is with the Internet and Data Science Lab (IDLab) of the Department of Information Technology (INTEC), Ghent University, B-9052 Ghent, Belgium e-mail: ivo.couckuyt@ugent.be.}%
	\thanks{J. Wauters, T. Lefebvre J. Degroote and G. Crevecoeur are member of core lab MIRO, Flanders Make, Belgium.}
	\thanks{I. Couckuyt is member of IMEC, Belgium.}
}

\markboth{}%
{}

\makeatletter
\let\@oldmaketitle\@maketitle
\renewcommand{\@maketitle}{\@oldmaketitle
	\vspace*{10pt}
	\hspace{6mm}
	\includegraphics[trim=0 6 0 6,clip,width=0.3\linewidth]{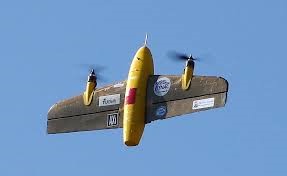}
	\includegraphics[trim=0 3 0 3,clip,width=0.3\linewidth]{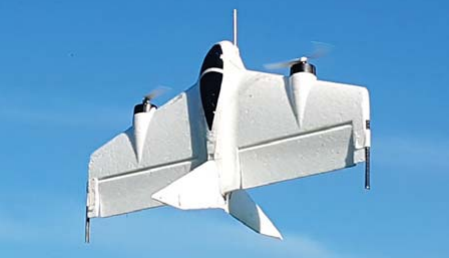}
	\includegraphics[trim=0 1 0 0,clip,width=0.3\linewidth]{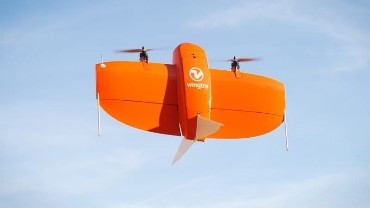}
	\vspace*{-0pt}}
\makeatother
\maketitle

\begin{abstract}
	
	In recent years, there has been a notable evolution in various multidisciplinary design methodologies for dynamic systems. Among these approaches, a noteworthy concept is that of concurrent conceptual and control design or co-design. This approach involves the tuning of feedforward and/or feedback control strategies in conjunction with the conceptual design of the dynamic system. The primary aim is to discover integrated solutions that surpass those attainable through a disjointed or decoupled approach. This concurrent design paradigm exhibits particular promise in the context of hybrid unmanned aerial systems (UASs), such as tail-sitters, where the objectives of versatility (driven by control considerations) and efficiency (influenced by conceptual design) often present conflicting demands. Nevertheless, a persistent challenge lies in the potential disparity between the theoretical models that underpin the design process and the real-world operational environment, the so-called reality gap. Such disparities can lead to suboptimal performance when the designed system is deployed in reality. To address this issue, this paper introduces DAIMYO, a novel design architecture that incorporates a high-fidelity environment, which emulates real-world conditions, into the procedure in pursuit of a `first-time-right' design. The outcome of this innovative approach is a design procedure that yields versatile and efficient UAS designs capable of withstanding the challenges posed by the reality gap.
	
\end{abstract}

\begin{IEEEkeywords}
Unmanned aerial vehicles, Multi-disciplinary design, Co-design, Bayesian optimization, Trajectory optimization, Differential flatness.
\end{IEEEkeywords}

\IEEEpeerreviewmaketitle

\section{Introduction}

\lettrine{L}{ast} decade has seen extensive exploration of hybrid \textsl{unmanned aerial system} (UAS) technologies combining fixed-wing and rotary-wing systems. Among these hybrids, the tail-sitter design stands out for its unique aerodynamic efficiency and the ability to transition between flight and hover by tilting its body ninety degrees \cite{saeed:2018,ducard:2021}. The tail-sitter is furthermore characterized by a \textsl{blended wing body} (BWB), a tailless design to increase its efficiency and typically two propellers. Stabilization is realized through a combined effort of differential thrust and placing control surfaces (elevons) in the propulsion wash (see the Figure at the top of the page with from left to right the Cyclone \cite{smeur:2019}, the X-VERT \cite{fuhrer:2019} and the Wingtra One \cite{olsson:2021}). Model-based system design plays a crucial role in advancing these systems, emphasizing the use of behavioural models for conceptual and control design, ultimately aiming for more cost-effective and innovative systems with enhanced functionalities, improved energy efficiency, and the capability to perform complex missions autonomously. 

Until recently, the development of tail-sitters has treated conceptual design, focused on the optimization of the aerodynamic efficiency of the BWB design \cite{theys:2016}, and control design, focused on feedback design (controller tuning) \cite{smeur:2019,fuhrer:2019} and to a limited extent feedforward design (trajectory optimization) \cite{banazadeh:2015, oosedo:2017, verling:2017, li:2020e}, separately. Within the field of \textsl{multi-disciplinary optimization} (MDO) the concept of concurrent conceptual and control design, or (control) \textsl{co-design}, connects these two aspects, expanding the design possibilities and yielding integrated solutions that cannot be achieved through a conventional isolated approach \cite{garcia:2019,herber:2018}. In recent years, there has been a growing focus on employing co-design methods in the realm of UAS development. This emphasis has shifted towards investigating architectural aspects, integrating both structural and propulsion aspects next to the conceptual and feedforward control design, mitigating computational expenses through the utilization of surrogate-supported techniques \cite{kaneko:2023, wauters:2022c, wauters:2023, balesdent:2022, morita:2020, matos:2022, lupp:2022, hendricks:2020, delbecq:2021}.

In general, two main architectures can be considered for the integration of conceptual design optimization and control design optimization: (1) a nested approach (also referred to as a `distributed' approach in the MDO research community \cite{martins:2013}), in which the optimal trajectory and/or feedback control parameters are determined for each iteration of the conceptual design optimization problem, and (2) a simultaneous approach (also referred to as a `monolithic' approach), in which the optimal conceptual design and optimal control design are determined together. 

The simultaneous approach directly extends on prevailing numerical schemes tailored to trajectory optimization \cite{ha:2018,whitman:2018,toussaint:2021}. Generally speaking the trajectory optimization problem is discretised (e.g. by parametrizing the trajectory, using a Direct Multiple Shooting discretization strategy, etc.). As such the problem is recast into a standard (constrained) \textsl{non-linear numerical optimization problem} (NLP). The resulting problem can then be addressed using appropriate NLP approaches. In the context of co-design, the optimization variables stemming from the trajectory parametrization are now extended with static design variables classifying it as a simultaneous co-design approach. Trajectory optimization problems are usually constrained, very often the NLP is solved using gradient based optimizers. However, the inclusion of design parameters in the optimization variables will then deteriorate the constraint sparsity of the NLP. This approach can be characterized by a reduced number of iterations and computational time, but has a larger and more complex design space and is therefore often confronted with convergence issues. 

\begin{figure}[t]
	\centering
	\includegraphics[trim=0 40 0 30,clip,width=0.8\linewidth]{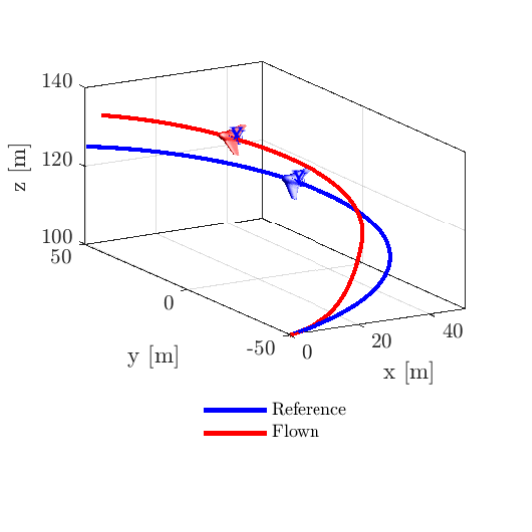}
	\caption{Visualization of the reality gap: closed loop flight trajectories of a low-fidelity (reference) and high fidelity (flown) model with identical reference signals.}	
	\label{fig:reality}
	\vspace{-5mm}
\end{figure}

On the other hand, the nested approach can be seen as a natural extension of conceptual design optimization in which conventional trajectory optimization is introduced as a constrained problem. As such, the optimization variables stemming from the trajectory parametrization are omitted from the outer loop NLP after they have been fixed in the inner loop. Therefore, it allows for separate definitions of static and dynamic objectives, as well as the use of solvers that are tailored to the specific characteristics of the problem. In doing so it avoids the convergence problems of the simultaneous approach. However, this approach is also subject to drawbacks, such as computational cost \cite{herber:2018}. This problem has recently been tackled through the inclusion of Bayesian optimization \cite{wauters:2022c, wauters:2023, lefebvre:2023}.

 However, a current shortcoming of the proposed frameworks is the sensitivity to \textsl{reality gap} (Fig.~\ref{fig:reality}), the offset between the physical system and the design model, due to the use of low fidelity models to enable trajectory optimization. To minimize this effect and to ensure that the resulting mission-specific designs are reliable in operation, variability and high-fidelity physics should be taken into account in the design procedure, as well as the computational cost that comes with it. This paper's main contribution lies in presenting DAIMYO, a first-time-right dynamic design architecture that closely aligns with needs and practices of practitioners, allowing for a smooth uptake by different communities. Additionally DAIMYO is applied to the design of a tail-sitter UAS, by leveraging on recent academic developments for tail-sitters such as flatness-based trajectory optimization \cite{tal:2023}, cascaded feedback linearization control \cite{ritz:2017} and fail-safe Bayesian optimization \cite{wauters:2024}. 

The remainder of this paper is organized as follows. In sec.~\ref{sec:codesign} of this paper, the dynamic design architecture is discussed during which the nested optimization loops are introduced. In sec.~\ref{sec:appl} the system modelling and control of a tail-sitter UAS is discussed. Finally, in sec.~\ref{sec:impl} the tail-sitter dynamics and control are introduced in the framework and applied on a dynamic design case to illustrate the working of the framework.

\section{First-time-right dynamic design architecture}\label{sec:codesign}

\begin{figure}[b]
	\centering
	\includegraphics[trim=0 0 0 0,clip,width=0.99\linewidth]{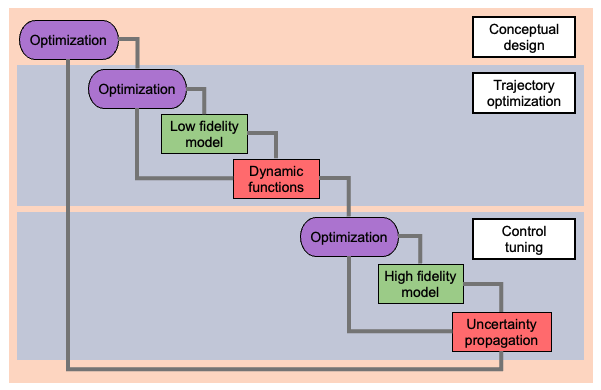}
	\caption{Simplified XDSM representation of DAIMYO.}	
	\label{fig:XDSM1}
\end{figure}

In this section we introduce DAIMYO, a first-time-right dynamic design framework. This multidisciplinary design architecture formulated as a nested/distributed co-design approach. The outer loop solves a conceptual design optimization problem of optimally tuned trajectories and controller settings. These trajectories and settings are realized in the inner loop which contains a sequential optimization of the flight paths and tuning of the controller. To pursue the first-time-right paradigm, a reality emulator is introduced by coupling a high-fidelity aerodynamic simulator with the rigid body dynamics of the system. This reality emulator is used to tune the controller and in doing so robustify the system against the reality gap. The framework is visualized through a simplified \textsl{Extended Design System Matrix} (XDSM \cite{martins:2013}) in Fig.~\ref{fig:XDSM1}. In what follows we move downstream through the waterfall structure to discuss the three optimization loops in more detail.

\subsection{Outer loop conceptual design}

The outer loop of our nested co-design framework solves a \text{conceptual design problem} (CDP), within which each design is evaluated for an optimal trajectory and optimally tuned control parameters. The design vector, $\vectorstyle{d}$, affects the trajectory optimization problem through the aerodynamic models. The performance of the UAS is evaluated in a reality emulator within which the UAS is flown closed-loop. The reality emulator contains a level of irreducible stochasticity to mimic varying operation conditions and is further elaborated on in section~\ref{sec:AVL}. Assessment of the performance will involve an \textsl{uncertainty propagation} (UP) step. We propose a multi-objective formulation of mean versus variance of the integrand cost rate $l_{\text{CDP}}$ over time domain $t\in[0,T]$. This corresponds to finding the Pareto front $\mathcal{P}$: the surface of non-dominated metrics on which we cannot improve on one without deteriorating on the other(s). 
\begin{equation}
	\begin{aligned}\label{eq:CDP}
		\matrixstyle{D}_*=\arg&\Par_{\substack{\matrixstyle{D}\subset\mathcal{D}}} ~ \{\mathbb{E}_{\vectorstyle{\epsilon}},\mathbb{V}_{\vectorstyle{\epsilon}}\}\left[\int\nolimits_{0}^Tl_{\text{\text{CDP}}}(\vectorstyle{\xi}(t),\vectorstyle{\upsilon}(t)|\matrixstyle{D}_*,\matrixstyle{\Xi}_*,\matrixstyle{K}_*)\text{d}t \right] \\
		&\begin{aligned}
			\text{ s.t. } \dot{\vectorstyle{\xi}}(t) &= \mathcal{F}_{\text{HF}}(\vectorstyle{\xi}(t),\vectorstyle{\upsilon}(t)|\vectorstyle{d},\vectorstyle{\epsilon}), && t\in[0,T]\\
			\vectorstyle{\upsilon}(t) &=\vectorstyle{\pi}(\vectorstyle{\xi}(t)|\vectorstyle{\xi}_*(t),\vectorstyle{\upsilon}_*(t),\vectorstyle{k}_*), && t\in[0,T]\\
		\end{aligned}
	\end{aligned}
\end{equation}
\noindent where $\vectorstyle{\xi}$ and $\vectorstyle{\upsilon}$ respectively denote the state and control vectors. Additionally, $\mathcal{F}_{\text{HF}}$ refers to a high-fidelity aerodynamic model corresponding to the reality emulator and $\vectorstyle{\epsilon}$ denotes the inherent stochasticity of high-fidelity environment. $\vectorstyle{\pi}$ represents the tuneable feedback control law as a function of, amongst others, the reference trajectory $\{\vectorstyle{\upsilon}_*,\vectorstyle{\xi}_*\}$ and the control parameters $\vectorstyle{k}_*$. Additionally, $\matrixstyle{D}\triangleq\{\vectorstyle{d}_i|i=1,...,n\}$ represents a set of design vectors and $\matrixstyle{D}_*$ corresponds to the Pareto optimal design vector set. Furthermore, $\mathcal{D}$ represents the design space of interest, often a compact subset of $\mathbb{R}^{d({\mathcal{D})}}$. For each design $\vectorstyle{d}$, an optimal trajectory $\{\vectorstyle{\upsilon}_*,\vectorstyle{\xi}_*\}$ exists, such that $\{\matrixstyle{\Upsilon}_*,\matrixstyle{\Xi}_*\}$ is the set of optimal trajectories corresponding to $\matrixstyle{D}_*$ and an optimal control parameter vector $\vectorstyle{k}_*$ exists, such that $\matrixstyle{K}_*$ is the set of optimal control parameters corresponding to $\matrixstyle{D}$. The solution of the dynamic design problem is therefore the Pareto set with corresponding trajectories $\{\matrixstyle{D}_*,\matrixstyle{\Upsilon}_*,\matrixstyle{\Xi}_*,\matrixstyle{K}_*\}$. 

\subsection{Inner loop trajectory optimization}\label{sec:to}

Dynamic system design, or co-design, differs from conventional MDO frameworks through its need to account for time-dependent behaviour. As such, critical to the design of highly dynamic systems is the need to determine the optimal behaviour of the system in time to realize a specific goal. To determine optimal flight trajectories we solve an open-loop \textsl{optimal control problem} (OCP)
\begin{equation}
	\label{eq:ocp}
	\begin{aligned}
		\vectorstyle{\xi}_*,\vectorstyle{\upsilon}_*=\arg&\min_{\vectorstyle{\xi},\vectorstyle{\upsilon}} \int_0^T l_{\text{\text{OCP}}}(\vectorstyle{\xi}(t),\vectorstyle{\upsilon}(t)|\vectorstyle{d}) \text{d}t \\
		&\begin{aligned}
			\text{ s.t. } \dot{\vectorstyle{\xi}}(t) &=\mathcal{F}_{\text{LF}}(\vectorstyle{\xi}(t),\vectorstyle{\upsilon}(t)|\vectorstyle{d}), && t\in[0,T]\\
			0 &= \vectorstyle{h}(\vectorstyle{\xi}(t),\vectorstyle{\upsilon}(t)), && t\in[0,T] \\
			0 &\geq \vectorstyle{g}(\vectorstyle{\xi}(t),\vectorstyle{\upsilon}(t)) , && t\in[0,T]
		\end{aligned}
	\end{aligned}
\end{equation}
\noindent were $\mathcal{F}_{\text{LF}}$ refers to a low-fidelity aerodynamic model. The functions $\vectorstyle{g}$ and $\vectorstyle{h}$ refer to path constraints. With slight abuse of notation we use the path equality constraint, $\vectorstyle{h}$, also to impose boundary conditions. The integrand cost rate, $l_{\text{OCP}}$, can impose any penalization deemed necessary. 

\subsection{Inner loop control tuning}\label{sec:tuning}

In pursuit of a first-time-right design, an effective controller needs to be realized that can account for the offset between the low fidelity models used for the trajectory optimization and a real-life setting. To realize this, a robust control strategy is proposed which is tuned in a reality emulator. This control strategy is further elaborated on section~\ref{sec:control} as it builds on the model of the tail-sitter. However, by definition of the framework, any tuneable control strategy can be implemented. 

To tune the controller, we propose a \textsl{feedback control problem} (FCP) in which the integrand cost rate, $l_{\text{FCP}}$, again contains any necessary form of penalization. As the emulator is stochastic, the integrand will be constructed out of terms derived from a UP approach (Fig.~\ref{fig:XDSM2}), which will be further elaborated on in section~\ref{sec:traj}.
\begin{equation}
	\label{eq:fcp}
	\begin{aligned}
		\vectorstyle{k}_*=\arg&\min_{\vectorstyle{k}} \mathbb{E}_{\vectorstyle{\epsilon}}\left[\int_0^T l_{\text{\text{FCP}}}(\vectorstyle{\xi}(t),\vectorstyle{\upsilon}(t)|\vectorstyle{d}) \text{d}t\right]\\
		&\begin{aligned}
			\text{ s.t. } \dot{\vectorstyle{\xi}}(t) &= \mathcal{F}_{\text{HF}}(\vectorstyle{\xi}(t),\vectorstyle{\upsilon}(t)|\vectorstyle{d},\vectorstyle{\epsilon}), && t\in[0,T]\\
			\vectorstyle{\upsilon}(t) &= \vectorstyle{\pi}(\vectorstyle{\xi}(t)|\vectorstyle{\xi}_*(t),\vectorstyle{\upsilon}_*(t),\vectorstyle{k}), && t\in[0,T]\\
		\end{aligned}
	\end{aligned}
\end{equation}

\section{Application case: tail-sitter UAS}\label{sec:appl}

\begin{figure}[b]
	\centering
	\includegraphics[trim=0 0 0 0,clip,width=0.8\linewidth]{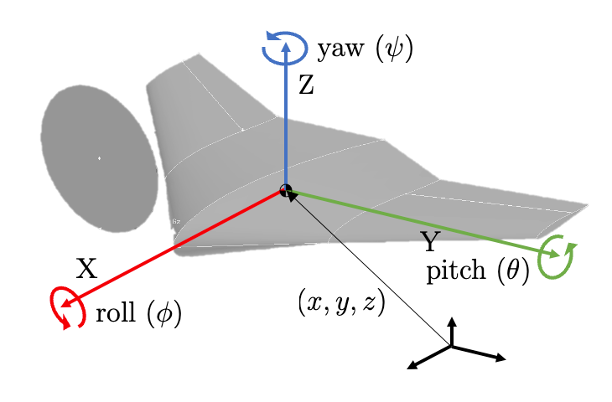}
	\caption{Body-fixed reference frame and angular motion terminology.}	
	\label{fig:reference}
	\vspace{-5mm}
\end{figure}

The DAIMYO framework was introduced in previous section in a generic form which allows the design of any dynamic system. In this paper, we are interested specifically in the design of tail-sitter UASs. Therefore, in this section we introduce the different aspects of the closed-loop dynamics of these systems.

In this work we use both a low-fidelity and high-fidelity description of the aerodynamic forces. The low-fidelity model is suited for trajectory optimization. The high-fidelity model is used to construct a simulation environment that allows high-fidelity estimates of the closed-loop performance but is too expensive for trajectory optimization purposes. 

\subsection{Rigid body dynamics}\label{sec:rbd}

We use Cartesian coordinates $x$, $y$ and $z$ gathered in the coordinate vector, $\vectorstyle{p}^{\text{global}}\in\mathbb{R}^3$, to represent the vehicle's position in a global inertial frame of reference. To represent the attitude of the system, we adopt the yaw-roll-pitch convention with yaw angle, $\psi$, roll angle, $\phi$, and, pitch angle, $\theta$. The angles are gathered in an angle vector $\vectorstyle{q}\in\mathbb{R}^3$ (Fig.~\ref{fig:reference}).

The flight dynamics of any rigid aircraft can be studied by means of the Newton-Euler equations of rigid body motion
\begin{subequations}
	\label{eq:rbd}
	\begin{align}
		m \dot{\vectorstyle{v}}^{\text{global}} + m \vectorstyle{g}^{\text{global}} &= \vectorstyle{f}^{\text{global}} \\
		\matrixstyle{I} \dot{\vectorstyle{\Omega}}^{\text{local}} + \vectorstyle{\Omega}^{\text{local}}\times\matrixstyle{I}\vectorstyle{\Omega}^{\text{local}} &=  \vectorstyle{\tau}^{\text{local}}
	\end{align}
\end{subequations}
Here $\vectorstyle{v}^{\text{global}}\in\mathbb{R}^3$ represents the linear velocity in the global frame of reference and $\vectorstyle{\Omega}^{\text{global}}\in\mathbb{R}^3$ represents the angular velocity in a body fixed local frame of reference. Likewise, $\vectorstyle{f}^{\text{global}}\in\mathbb{R}^3$ and $\vectorstyle{\tau}^{\text{local}}\in\mathbb{R}^3$ represent the total force and torque acting on or about the centre of mass, respectively expressed in the global and local frame of references.

The angular velocity $\vectorstyle{\Omega}^{\text{local}}$ is related to the time derivative of the angles, $\dot{\vectorstyle{q}}$ through the Jacobian matrix $\matrixstyle{J}(\vectorstyle{q})$. The $\vectorstyle{v}^{\text{global}}$ and time derivative of the position vector, $\dot{\vectorstyle{p}}^{\text{global}}$, are equal. The state of the system is stored in a vector $\vectorstyle{\xi}\in\mathbb{R}^{12}$ where $\vectorstyle{\xi}=(\vectorstyle{p}^{\text{global}},\vectorstyle{q},\vectorstyle{v}^{\text{global}},\vectorstyle{\Omega}^{\text{local}})$. 

Clearly, in the context of UASs, the force, $\vectorstyle{f}$, and torque, $\vectorstyle{\tau}$, are generated by the aerodynamic surfaces and the propulsion system. We examine an aircraft characterized by four control inputs, namely two propellers, generating thrust forces, $\text{T}_1$ and $\text{T}_2$, and, two control surfaces, $\delta_1$ and $\delta_2$. These four control inputs are stored in the input, $\vectorstyle{\upsilon}\in\mathbb{R}^4$ where ${\vectorstyle{\upsilon}} = (\text{T}_1,\text{T}_2, \delta_1 , \delta_2)$.

Further we utilise both a low- and high-fidelity description of the aerodynamic model. These are discussed next. We may note here that the aerodynamic characteristics of the system depend highly on the conceptual design, parametrized by design vector $\vectorstyle{d}$. This will not be made explicit in this section but is of critical importance.

\subsection{High-fidelity aerodynamic model}\label{sec:AVL}

As motivated in the introduction, we pursue an optimization framework with as little of a reality gap that can be evaluated with the computational budget available. To that end, a numerical solver is constructed that evaluates the aerodynamic forces $\vectorstyle{f}^{\text{local}}$ and $\vectorstyle{\tau}^{\text{local}}$ as a function of the system's state. For that matter, any computational fluid dynamics methodology can be used. In this work we rely on the aerodynamic forces predicted by Drela's Athena Vortice Lattice (AVL) \cite{drela:2017}.

AVL assumes on the one hand that the flow is inviscid and irrotational and on the other hand that the aircraft can be represented as a zero-thickness sheet of a finite number of panels. The combination of which basically implies that AVL can be seen as a discretized Prandtl lifting line method. The assumption of an inviscid and irrotational flow comes with a velocity potential, such that velocity field can be understood as the gradient of a scalar. Since under these conditions the superposition principle holds, namely that each sum of solutions is again a solution, a perturbation velocity potential can be introduced in each panel as an equivalent to the horse shoe vortices of the lifting line theory. The impact of each perturbation potential on each other potential can be calculated by means of Biot-Savart's law. Determining the flow around the aircraft now boils down to determining the strength of these individual potentials. This can be realized by enforcing the flow to be tangential to the surface. The results is a linear set of equations that can be efficiently solved. By means of Kutta-Joukowski's theorem, the force coefficients become readily available \cite{drela:2014}.

AVL provides the $\pi$-theory force coefficients, $C_X$, $C_Y$ and $C_Z$, and moment coefficients, $C_{\phi}$, $C_{\theta}$ and $C_{\psi}$, as a function of the flight direction, $\alpha$ and $\beta$, control deflections, $\delta_1$ and $\delta_2$, and rotational speed in the body frame $\vectorstyle{\Omega}^{\text{local}}$. Since AVL relies on the potential flow method, the force coefficients are Reynolds number independent. The forces and moments can now be obtained 
\begin{subequations}
	\label{eq:AVL}
	\begin{align}
		\vectorstyle{f}^\text{local}&= \begin{pmatrix}
			\text{T}_1+\text{T}_2 - C_{\text{X}} \sfrac{1}{2}\rho ||\vectorstyle{v}||^2 S \\
			C_{\text{Y}} \sfrac{1}{2}\rho ||\vectorstyle{v}||^2 S \\
			- C_{\text{Z}} \sfrac{1}{2}\rho ||\vectorstyle{v}||^2 S
		\end{pmatrix} \\
		\vectorstyle{\tau}^\text{local} &= \begin{pmatrix}
			C_{\phi} \sfrac{1}{2}\rho ||\vectorstyle{v}||^2 S c \\
			C_{\theta} \sfrac{1}{2}\rho ||\vectorstyle{v}||^2 S b \\
			C_{\psi} \sfrac{1}{2}\rho ||\vectorstyle{v}||^2 S c + l^{\text{T}}_y (\text{T}_2 -\text{T}_1)
		\end{pmatrix}
	\end{align}
\end{subequations}
\noindent $\rho$ represents the air density, $S$ the aircraft's surface, $b$ the span and $c$ the mean aerodynamic chord and $l^{\text{T}}_y$ the $y$-axis thrust moment arm.

While the computational cost of AVL is fairly low, directly coupling it do the rigid body dynamics to be solved is still computationally unfeasible. Therefore, cheap-to-sample surrogates of the force coefficients are constructed for each design as a function of the states. The surrogate model used is the Gaussian process interpolator, which is elaborated on in Appendix \ref{sec:BO}. A design of experiments containing a hundred states is used to construct the force coefficient surrogates.

\subsection{Low-fidelity aerodynamic model}\label{sec:forces}

The use of a high-fidelity model as discussed in the previous section is inconvenient when addressing problems such as trajectory generation. The main reason is, even with the computationally efficient surrogates in place, the dynamic equations are in a form that does not allow any analytical treatment. Hence, to obtain a sufficiently simple analytical model, a low-fidelity model, we adopt $\phi$-theory \cite{lustosa:2018}. The resulting model clearly presents an over-simplification, however it does capture the principal dynamic behavioural tendencies of the system, which is more than sufficient for its intended purpose. Specifically, the properties of the model make it ideally suited for trajectory optimization and control design.

A simplified representation is given below in which $l^{\text{T}}_y$ denotes the distance of the propellers from the symmetry plane and $K_{\text{L}}$, $K_{\text{D}}$, $K_{\phi}$, $K_{\theta}$ and $K_{\psi}$ respectively the lift, drag, roll, pitch and yaw coefficients.
\begin{subequations}
	\label{eq:forces}
	\begin{align}
		\vectorstyle{f}^\text{local}&= \begin{pmatrix}
			\text{T}_1+\text{T}_2 - K_{\text{D}} \vectorstyle{v}_x^\text{local} \|\vectorstyle{v}\| \\
			0 \\
			- K_{\text{L}} \vectorstyle{v}_z^\text{local} \|\vectorstyle{v}\|
		\end{pmatrix} \\
		\vectorstyle{\tau}^\text{local} &= \begin{pmatrix}
			K_{\phi} \vectorstyle{v}_x^\text{local} \|\vectorstyle{v}\| (\delta_2 - \delta_1) \\
			K_{\theta} \vectorstyle{v}_x^\text{local} \|\vectorstyle{v}\| (\delta_1 + \delta_2) \\
			K_{\psi} \vectorstyle{v}_x^\text{local} \|\vectorstyle{v}\| (\delta_2 - \delta_1) + l^{\text{T}}_y (\text{T}_2 -\text{T}_1)
		\end{pmatrix}
	\end{align}
\end{subequations}

\begin{figure*}[t]
	\centering
	\includegraphics[trim=0 0 0 0,clip,width=0.99\linewidth]{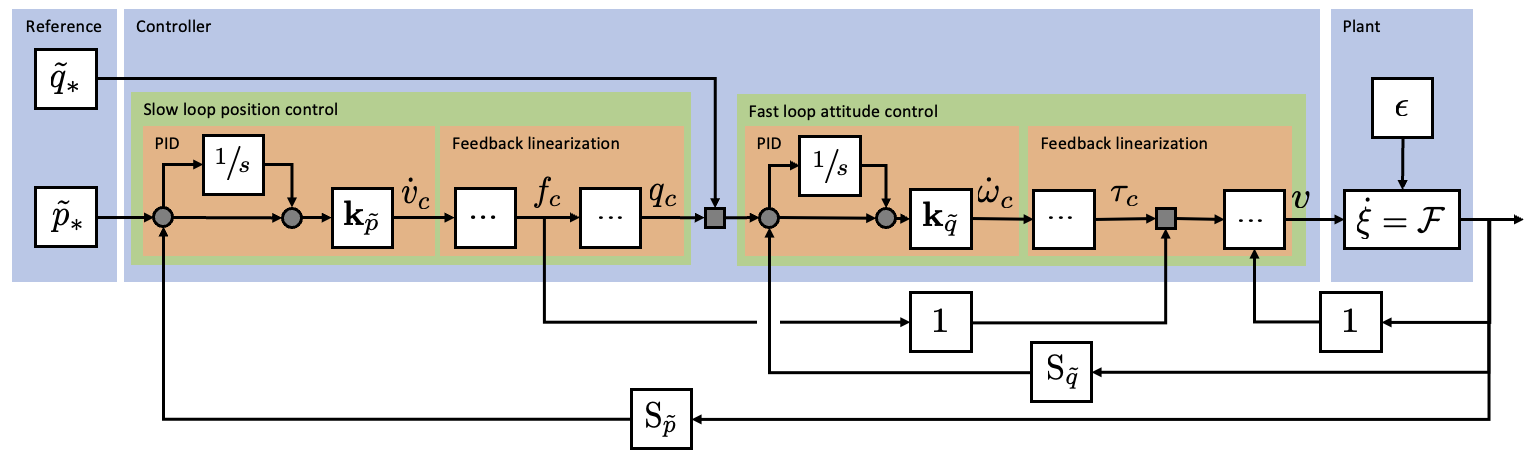}
	\caption{Control diagram of the novel cascaded dynamic feedback linearization controller.}	
	\label{fig:control}
	\vspace{-5mm}
\end{figure*}

To determine the coefficients of the low-fidelity model, again we use AVL. The force and moment coefficients introduced above rely on Lustosa's et al. \cite{lustosa:2018} $\phi$-theory, but can be linked to the more conventional coefficients derived using Buckingham $\pi$-theory that are provided by AVL: $K_{\text{L}}$ denotes the change in lift due to velocity, which can be approximated by 
\begin{equation}
	K_{\text{L}}\approx\sfrac{1}{2}\rho S\partial_\alpha C_{\text{L}}
\end{equation}
In a similar manner $K_{\text{D}}$ can be determined. However, Tal et al. argue that at low velocities (UAS operational range), this term approximates zero, as it is dominated by viscous effects \cite{tal:2023}. Since AVL relies on the potential flow theory, it will be set to zero. The moment coefficients $K_{\phi}$, $K_{\theta}$ and $K_{\psi}$ respectively correspond to the $\phi$-theory defined roll moment due to velocity and control deflection, which can be linked to $\pi$-theory as
\begin{subequations}
	\begin{align}
		K_{\phi}&\approx	\sfrac{1}{2}\rho S l^\delta_y \partial^2_{v,\delta} C_{\text{L}} \\
		K_{\theta}&\approx\sfrac{1}{2}\rho S l^\delta_x \partial^2_{v,\delta} C_{\text{L}} \\
		K_{\psi}&\approx	\sfrac{1}{2}\rho S l^\delta_y \partial^2_{v,\delta} C_{\text{D}}
	\end{align}
\end{subequations}
with $l^\delta_y$ the $y$-axis control surface moment arm and $l^\delta_x$ the $x$-axis moment arm. Tal et al. approximate $\partial^2_{v,\delta} C_{\text{L}}$ by $f_{\text{con}}  \partial_v C_{\text{L}}$ with $f_{\text{con}} $ the control fraction of the wing.

It is important to emphasize that the aerodynamic force in this context lacks any lateral component. This characteristic arises primarily from the fact that the majority of tail-sitter designs discussed in the existing literature do not feature a fuselage and are devoid of a tail section or any vertical aerodynamic elements. As shown by Tal \cite{tal:2023}, it follows that the tailsitter is differentially flat, an advantageous property for system control design. A system exhibits differential flatness when the states and inputs, the differentially dependent system variables, can be written in terms of a set of variables and their time derivatives, referred to as the flat output. The tail-sitter considered in this work is a flat with coordinates $\vectorstyle{\sigma}= (x,y,z,\psi)\in\mathbb{R}^4$. The differential flatness of the tailsitter system is further discussed in Appendix \ref{sec:flat}.

\subsection{Feedforward control design}\label{sec:forward}

Our trajectory optimization framework makes use of the fact that, as a direct result of the differential flatness, any feasible state-input trajectory, i.e. satisfying Eq.~\ref{eq:rbd} and making use of Eq.~\ref{eq:forces}, can be related to an equivalent signal in the flat output, $\vectorstyle{\sigma}(t)$. To parametrize the flat trajectory, $\vectorstyle{\sigma}(t)$, we utilize the B-spline framework. The flat trajectory is defined as a weighted sum of B-splines basis functions, such that $\vectorstyle{\sigma}(t|\vectorstyle{c}) \triangleq \sum\nolimits_{i=0}^n B_{i,d}(t) \cdot\vectorstyle{c}_i$ with $\vectorstyle{c}_i$ the weighting coefficients that are to be determined. B-splines of order $d+1$ are function bases constructed out of piecewise polynomials $B_{i,d}$ with order $d$ that meet in a set of non-decreasing time instances $\mathcal{T}_{\text{splines}} ={t_0,t_1,\dots,t_m}$, referred to as the knot vector. When $t_0=0$ and $t_m=T$ the B-spline covers the time interval $[0,T]$ \cite{stoical:2016}. 

To evaluate the objective function we utilize a trapezoidal integration rule over a sufficiently dense time grid, $\mathcal{T}_{\text{int}}$. The same grid is used to evaluate the inequality path constraints. We only consider boundary conditions as equality path constraints. These can be evaluated without further approximation. These steps transcribe (\ref{eq:ocp}) into a non-linear program that can be solved using standard gradient based optimizers. 

\begin{equation}
	\begin{aligned}
		\vectorstyle{c}_* = \arg&\min_{\vectorstyle{c}} \sum\nolimits_{t\in\mathcal{T}_{\text{int}}} l_{\text{OCP}}(t|\vectorstyle{d},\vectorstyle{c}) \\
		&\text{s.t. }0\geq\vectorstyle{g}(t|\vectorstyle{d},\vectorstyle{c}) \\	
	\end{aligned}
\end{equation}

The optimal state and control signals are now directly obtained from $\vectorstyle{\xi}_* = \Xi[\vectorstyle{\sigma}(\vectorstyle{c}_*)]$ and $\vectorstyle{\upsilon}_* = \Upsilon[\vectorstyle{\sigma}(\vectorstyle{c}_*)]$ as described in Appendix \ref{sec:flat}.

\subsection{Feedback control design}\label{sec:control}

To simulate the closed-loop performance of the system we introduce a control strategy based on the work of Ritz et al. \cite{ritz:2017}. The control design relies heavily on the simplified dynamics from section \ref{sec:forces}. The controller can be understood as a cascaded \textsl{static feedback linearization}. The cascading refers to the exploitation of the so-called time-scale separation principle of attitude and position, which allows the use of two nested control structures operating at different time-scales.

First a desired force, $\vectorstyle{f}_c$, is calculated according the feedback linearization principle and disregarding the fact that we cannot generate just any force. We add an integration term for robustness where $\dot{\vectorstyle{\gamma}}=\vectorstyle{\xi}-\vectorstyle{\xi}_*$. This part determines the outer-loop. 
\begin{subequations}
	\label{eq:control_f}
	\begin{align}
		\vectorstyle{f}_c(\vectorstyle{\xi},\vectorstyle{\xi}_*,\dot{\vectorstyle{\xi}}_*,\vectorstyle{\gamma}) = &\: m \dot{\vectorstyle{v}}_c(\vectorstyle{\xi},\vectorstyle{\xi}_*,\dot{\vectorstyle{\xi}}_*,\vectorstyle{\gamma}) + m \vectorstyle{g} \\
		\dot{\vectorstyle{v}}_c(\vectorstyle{\xi},\vectorstyle{\xi}_*,\dot{\vectorstyle{\xi}}_*,\vectorstyle{\gamma}) = &\: \dot{\vectorstyle{v}}_* - \vectorstyle{k}_p^2 (\vectorstyle{v} - \vectorstyle{v}_*) \nonumber\\
		&\:- \vectorstyle{k}^1_p (\vectorstyle{p} - \vectorstyle{p}_*) - \vectorstyle{k}^0_p \vectorstyle{\gamma}_p
	\end{align}
\end{subequations}

Then, we calculate an instantaneous desired attitude, $\matrixstyle{R}_c\in\mathfrak{SO}(3)$. The feasible aerodynamic force, $\vectorstyle{f}$, should be matched with the desired force, $\vectorstyle{f}_c$. Using Ritz et al.'s additional condition $\psi_* = \arctan \sfrac{\dot{y}_*}{\dot{x}_*}$ the target attitude, $\matrixstyle{R}_c$, is calculated. The corresponding desired roll and pitch follow from the flatness expressions. We refer to Appendix \ref{sec:flat} for details. Once the control attitude, $\matrixstyle{R}_c$, is available, a desired moment, $\vectorstyle{\tau}_c$, is calculated again according the static feedback linearization principle.
\begin{equation}
	\label{eq:control_tau}
	\vectorstyle{\tau}_c(\vectorstyle{\xi},\vectorstyle{\xi}_*,\dot{\vectorstyle{\xi}}_*,\vectorstyle{\gamma}) = \matrixstyle{I} \dot{\vectorstyle{\Omega}}_c(\vectorstyle{\xi},\vectorstyle{\xi}_*,\dot{\vectorstyle{\xi}}_*,\vectorstyle{\gamma}) + {\vectorstyle{\Omega}}\times\matrixstyle{I}{\vectorstyle{\Omega}} 
\end{equation}

Inspired by Ritz et al.'s look-up table assisted non-linear \textsl{model predictive control} (MPC) strategy, the following empirically derived control law is proposed to calculate the control velocity, $\dot{\vectorstyle{\Omega}}_c$
\begin{equation}
	\dot{\vectorstyle{\Omega}}_c(\vectorstyle{\xi},\vectorstyle{\xi}_*,\dot{\vectorstyle{\xi}}_*,\vectorstyle{\gamma}) =\vectorstyle{k}^2_q ({\vectorstyle{\Omega}}-\vectorstyle{\Omega}_*) + \vectorstyle{k}_q^1 \Delta \vectorstyle{q} + \vectorstyle{k}_q^0\vectorstyle{\gamma}_q
\end{equation}

By means of the flat inverse dynamics, described in \ref{sec:flat}, the control inputs $\text{T}_1$, $\text{T}_2$, $\delta_1$ and $\delta_2$ can be obtained from the control force and torque, $\vectorstyle{f}_c$ and $\vectorstyle{\tau}_c$ respectively. In accordance with the inherent behaviour of the system, the control should pursue obtaining the desired attitude first, effectively aligning the control force, $\vectorstyle{f}_c$, with the feasible force plane. This can be realized by ensuring that the attitude controller's time constant is a magnitude smaller than the position controller's one. This is an important consideration when tuning the gain matrices $\{\vectorstyle{k}_j^i|j=p,q;i=0,1,2\}$. The tuning of the control parameters is discussed in \ref{sec:tuning}. 

A control diagram of the architecture is presented in Fig.~\ref{fig:control}. One observes that attitude and position are tracked using a cascade of feedback linearizations in conjunction with PID control strategies. To visualize this $\matrixstyle{S}$ is introduced as a selection of part of the state. Additionally $\tilde{\cdot}$ denotes the quantity of interest and its derivatives. Finally, $\sfrac{1}{s}$ is borrowed from the conventional frequency domain terminology to denote the integration of the state offset over time.

\section{Implementation and validation}\label{sec:impl}

A dynamic design problem of a tail-sitter is addressed with the incentive of validating DAIMYO and discussing the implementation. The control -- both forward and backward -- and conceptual design are optimized towards the optimal performance of the abstraction of real-life application in the form of a discrete set of missions. The mission profiles under consideration, the tail-sitter's design parametrization, the discussion of the implementation of the design and control optimization, the outcomes derived from our concurrent conceptual and control design optimization architecture, and an ablation study to emphasize the importance of dynamic design are discussed next.

\subsection{Mission specifications}\label{sec:mission} 

We define four lower level trajectory optimization problems that aim to capture the essential functionalities that should be displayed by the UAS: a flight transition from vertical to horizontal flight and vice versa, mimicking the take-off and landing, cruise flight and a turn manoeuvre. The cruise flight and a turn manoeuvre are defined from the perspective that once their corresponding objectives are optimized, the tail-sitter will be enabled to perform Dubins path motivated planning in an efficient manner. The key idea behind Dubins paths is to find the shortest path between two points in a three-dimensional space while considering the vehicle's limited turning radius or turning rate. The optimal path can be composed of any combination of a straight line and a circular arc \cite{dubins:1957,owen:2014b}. Finding the optimal Dubins path involves solving a mathematical problem to determine the sequence and lengths of these segments. 

The problem is defined as a two-step optimization. In the first step, the boundary conditions of the trajectory optimization problem are obtained. This is realized by determining the time span, $T$, and velocities, $\vectorstyle{bc}\triangleq\{\vectorstyle{v}(0), \vectorstyle{v}(T), \vectorstyle{\Omega}(0), \vectorstyle{\Omega}(T)\}$, such that for a fixed initial and final position, $\vectorstyle{p}(0)$ and $\vectorstyle{p}(T)$, and orientation, $\vectorstyle{q}(0)$ and $\vectorstyle{q}(T)$, a third order polynomial $P_3$ for $\vectorstyle{\sigma}(t)$ can be defined for which the thrust $\text{T}_1$ and $\text{T}_2$ is minimal.
\begin{equation}
	\label{eq:bc}
	\begin{aligned}
		\vectorstyle{bc}_*=\arg&\min_{\vectorstyle{bc}} \int_0^T (\text{T}_1+\text{T}_2)\text{d}t \\
		&\begin{aligned}
			\text{ s.t. } \vectorstyle{\sigma}(t) &=P_3(\vectorstyle{bc}), && t\in[0,T]\\
			0 &\geq \vectorstyle{g}(\vectorstyle{\xi}(t),\vectorstyle{\upsilon}(t)) , && t\in[0,T] \\	
		\end{aligned}
	\end{aligned}
\end{equation}

\noindent The inequality constraint $\vectorstyle{g}(\vectorstyle{\xi}(t),\vectorstyle{\upsilon}(t))$ ensures that the maximal enforceable control signals $\{\text{T}_{\text{min}},\text{T}_{\text{max}},\delta_{\text{min}},\delta_{\text{max}}\}$ are not exceeded. The second step involves the flatness-based trajectory optimization (cf. sec.~\ref{sec:forward}) by means of the optimized boundary conditions. 

\subsection{Conceptual design parameterization} 

The BWB can be conceptually divided into two distinct wing sections, each characterized by its airfoil type, chord length ($c$), leading edge position in the body frame ($x$, $y$, and $z$), and orientation relative to the chord on the symmetry plane at three critical locations: the wingtip (denoted as subscript tip), the transition point from the wing to the fuselage (subscript fus), and the aircraft's symmetry plane (subscript sym). The symmetry plane's leading edge serves as the reference point for these three sections. It's important to note that in this study, we maintain a fixed section profile corresponding to the CAL4014l flying wing reflexed airfoil \cite{williamson:2012}. The aircraft's propulsion system is positioned at the leading edge of the outer wing's center, while the control surfaces are defined to span the trailing edge of the outer wing, with their leverage point at the center such that the lever arms, $l^{\text{T}}_x$, $l^{\text{T}}_y$, $l^\delta_x$ and $l^\delta_y$, are fully defined. The size of these control surfaces can be adjusted by specifying the fraction of the chord at the tip where the control surface starts ($f_{\text{con}} $). Using these $d(\vectorstyle{d})=12$ design variables a comprehensive definition of the aircraft's configuration is readily available \cite{sobester:2014} and permits the calculation of all other parameters required to evaluate the aerodynamics and rigid body dynamics (Table~\ref{table:variables}).

\begin{figure}[t]
	\centering
	\includegraphics[trim=0 0 0 0,clip,width=0.8\linewidth]{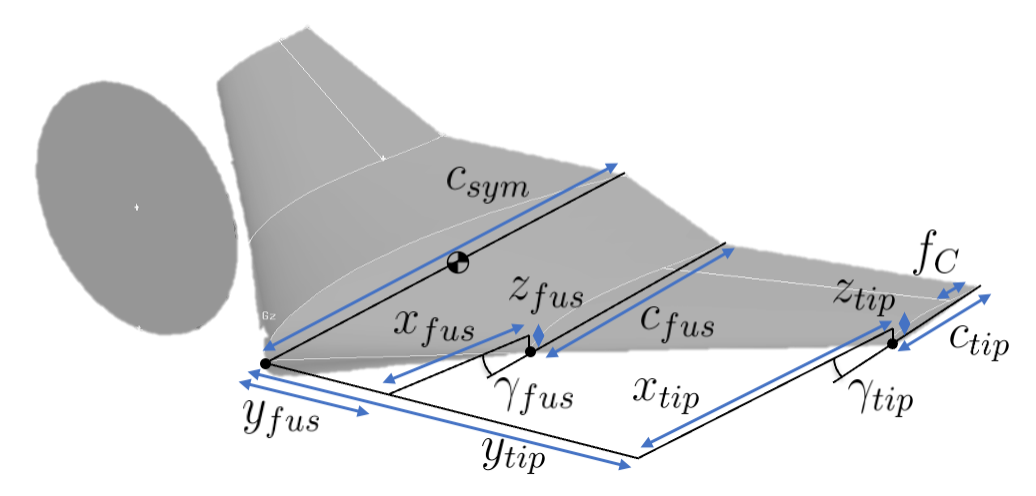}
	\caption{Parameterization of the tail-sitter.}	
	\label{fig:ux5}
	\vspace{-5mm}
\end{figure}

\begin{landscape}
	\begin{figure}[t]
		\centering
		\includegraphics[trim=0 0 0 0,clip,width=0.99\linewidth]{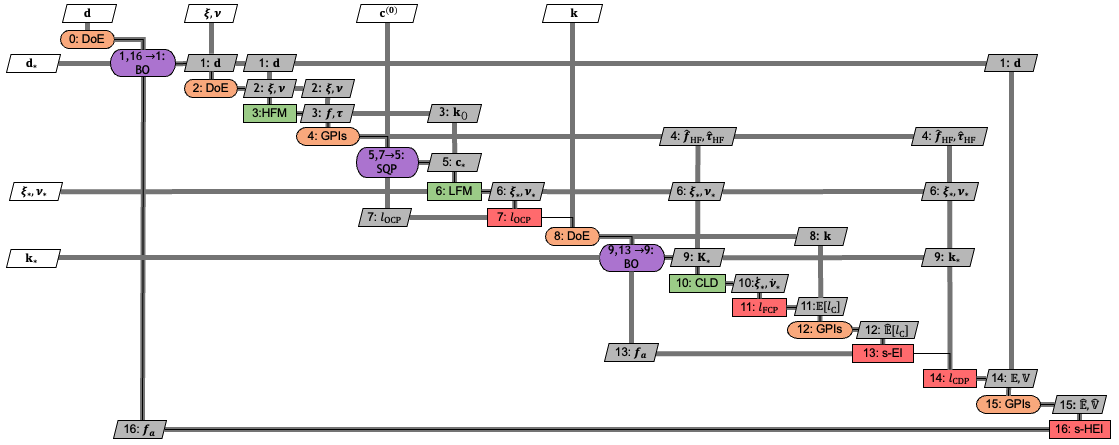}
		\caption{Extended Design System Matrix (XDSM \cite{martins:2013}) of the first-time-right dynamic design framework DAIMYO. Grey lines represent the data flow, black lines the process flow, horizontal lines are output, vertical lines are input. We use a $\wedge$ above a quantity of interested that has been modelled with a GPI.}	
		\label{fig:XDSM2}
	\end{figure}
\end{landscape}

\begin{table*}[p]
	\centering
	\renewcommand{\arraystretch}{1.4}
	\begin{tabular}{ r l | r r | r r r r | r r}
		\multicolumn{2}{c |}{design} & lower  & upper  &  &  &  &  & ablation   & ablation \\
		\multicolumn{2}{c |}{variables} &  bound &  bound & baseline & design 1 & design 2 & design 3 & design 1 &  design 2 \\
		\hline
		$c_{\text{sym}}$ & [$\si{m}$] & $0.4000$ & $0.8000$ & $0.6000$ & $0.6625$ & $0.7080$ & $0.6375$ & $0.4545$ & $0.4084$ \\
		$c_{\text{fus}}$ & [$\si{m}$] & $0.2000$ & $0.4000$ & $0.3000$ & $0.2791$ & $0.3368$ & $0.2018$ & $0.2000$ & $0.2000$ \\
		$c_{\text{tip}}$ & [$\si{m}$] & $0.1000$ & $0.2000$ & $0.1500$ & $0.1301$ & $0.1057$ & $0.1097$ & $0.1000$ & $0.1000$ \\
		$x_{\text{fus}}$ & [$\si{m}$] & $0.1000$ & $0.3000$ & $0.2000$ & $0.2022$ & $0.1659$ & $0.2338$ & $0.1000$ & $0.1000$ \\
		$x_{\text{tip}}$ & [$\si{m}$] & $0.4000$ & $0.6000$ & $0.5000$ & $0.4027$ & $0.4952$ & $0.4576$ & $0.4000$ & $0.4000$ \\
		$y_{\text{fus}}$ & [$\si{m}$] & $0.1000$ & $0.3000$ & $0.2000$ & $0.1721$ & $0.1453$ & $0.2105$ & $0.1000$ & $0.1001$ \\
		$y_{\text{tip}}$ & [$\si{m}$] & $0.4000$ & $0.6000$ & $0.5000$ & $0.5809$ & $0.5217$ & $0.5733$ & $0.5183$ & $0.4000$ \\
		$z_{\text{fus}}$ & [$\si{m}$] & $-0.0100$ & $0.0100$ & $0.0000$ & $0.0026$ & $-0.0097$ & $0.0098$ & $-0.0100$ & $0.0031$ \\
		$z_{\text{tip}}$ & [$\si{m}$] & $-0.0100$ & $0.0100$ & $0.0000$ & $0.0100$ & $-0.0085$ & $-0.0094$ & $0.0015$ & $-0.0100$ \\
		$\gamma_{\text{fus}}$ & [$^\circ$] & $-6.0000$ & $4.0000$ & $-1.0000$ & $-3.8612$ & $1.1102$ & $-2.5758$ & $-2.4010$ & $-1.9364$ \\
		$\gamma_{\text{tip}}$ & [$^\circ$] & $-12.0000$ & $-2.0000$ & $-7.0000$ & $-2.7843$ & $-8.2476$ & $-2.0626$ & $-2.0001$ & $-2.9237$ \\
		$f_{\text{con}} $ & [$-$] & $0.5000$ & $0.7500$ & $0.6250$ & $0.5145$ & $0.7326$ & $0.7087$ & $0.5009$ & $0.7250$ \\
  		$\vectorstyle{k}^q_p$ & [$-$] & $0.0000$ & $2000.0000$ & $1508.0752$ & $1194.7066$ & $1657.4292$ & $813.9534$ & $-$ & $-$ \\
		$\vectorstyle{k}^q_d$ & [$-$] & $0.0000$ & $40.0000$ & $12.0000$ & $5.0920$ & $3.1903$ & $9.3587$ & $-$ & $-$ \\
		$\vectorstyle{k}^q_i$ & [$-$] & $0.0000$ & $2000.0000$ & $467.7640$ & $1908.3785$ & $172.2978$ & $1937.5000$ & $-$ & $-$ \\
		$\vectorstyle{k}^p_p$ & [$-$]  & $0.0000$ & $0.8000$ & $0.0091$ & $0.5321$ & $0.0668$ & $0.0500$ & $-$ & $-$ \\
		$\vectorstyle{k}^q_d$ & [$-$] & $0.0000$ & $4.0000$ & $0.1817$ & $0.7364$ & $1.1824$ & $0.9750$ & $-$ & $-$ \\
		$\vectorstyle{k}^q_i$ & [$-$] & $0.0000$ & $0.8000$ & $0.0000$ & $0.0000$ & $0.0000$ & $0.0000$ & $-$ & $-$ \\
		\hline
		$m$ & [$\si{kg}$] & $1.3893$ & $4.3368$ & $2.3764$ & $2.3869$ & $2.4655$ & $2.1820$ & $1.5010$ & $1.3962$ \\
		$x_{\text{cog}}$ & [$\si{m}$] & $0.2000$ & $0.4000$ & $0.3000$ & $0.3312$ & $0.3540$ & $0.3188$ & $0.2272$ & $0.2042$ \\
		$y_{\text{cog}}$ & [$\si{m}$] & $0.0000$ & $0.0000$ & $0.0000$ & $0.0000$ & $0.0000$ & $0.0000$ & $0.0000$ & $0.0000$ \\
		$z_{\text{cog}}$ & [$\si{m}$] & $0.0000$ & $0.0000$ & $0.0000$ & $0.0000$ & $0.0000$ & $0.0000$ & $0.0000$ & $0.0000$ \\
		$l^{\text{T}}_x$ & [$\si{m}$] & $-0.0500$ & $-0.0500$ & $-0.0500$ & $0.0288$ & $0.0234$ & $-0.0269$ & $-0.0228$ & $-0.0458$ \\
		$l^{\text{T}}_y$ & [$\si{m}$] & $0.2500$ & $0.4500$ & $0.3500$ & $0.3765$ & $0.3335$ & $0.3919$ & $0.3091$ & $0.2500$ \\
		$l^{\delta}_x$ & [$\si{m}$] & $-0.2000$ & $-0.3500$ & $-0.2750$ & $-0.1758$ & $-0.1978$ & $-0.1827$ & $-0.1728$ & $-0.1958$ \\
		$l^{\delta}_y$ & [$\si{m}$] & $0.2500$ & $0.4500$ & $0.3500$ & $0.3765$ & $0.3335$ & $0.3919$ & $0.3091$ & $0.2500$ \\
		$\text{I}_{\text{xx}}$ & [$\si{kgm^2}$] & $0.0048$ & $0.0874$ & $0.0246$ & $0.0283$ & $0.0222$ & $0.0185$ & $0.0098$ & $0.0047$ \\
		$\text{I}_{\text{yy}}$ & [$\si{kgm^2}$] & $0.0011$ & $0.0113$ & $0.0039$ & $0.0018$ & $0.0031$ & $0.0019$ & $0.0012$ & $0.0011$ \\
		$\text{I}_{\text{zz}}$ & [$\si{kgm^2}$] & $0.0058$ & $0.0982$ & $0.0284$ & $0.0301$ & $0.0252$ & $0.0204$ & $0.0110$ & $0.0058$ \\
		$\text{I}_{\text{xy}}$ & [$\si{kgm^2}$] & $0.0000$ & $0.0000$ & $0.0000$ & $0.0000$ & $0.0000$ & $0.0000$ & $0.0000$ & $0.0000$ \\
		$\text{I}_{\text{xz}}$ & [$\si{kgm^2}$] & $0.0000$ & $0.0000$ & $0.0000$ & $0.0000$ & $0.0000$ & $0.0000$ & $0.0000$ & $0.0000$ \\
		$\text{I}_{\text{yz}}$ & [$\si{kgm^2}$] & $0.0000$ & $0.0000$ & $0.0000$ & $0.0000$ & $0.0000$ & $0.0000$ & $0.0000$ & $0.0000$ \\
		$S$ & [$\si{m^2}$] & $0.1500$ & $0.5400$ & $0.3150$ & $0.3294$ & $0.3183$ & $0.2897$ & $0.1909$ & $0.1509$ \\
		$b$ & [$\si{m}$] & $0.8000$ & $1.2000$ & $1.0000$ & $1.1618$ & $1.0434$ & $1.1467$ & $1.0366$ & $0.8000$ \\
		$c$ & [$\si{m}$] & $0.1875$ & $0.4500$ & $0.3150$ & $0.2835$ & $0.3051$ & $0.2527$ & $0.1842$ & $0.1886$ \\
		$C_{\text{D,0}}$ & [$-$] & $0.0003$ & $0.0006$ & $0.0002$ & $0.0001$ & $0.0004$ & $0.0001$ & $0.0000$ & $0.0000$ \\
		$k$ & [$-$] & $0.0835$ & $0.1147$ & $0.0978$ & $0.0795$ & $0.0942$ & $0.0710$ & $0.0563$ & $0.0772$ \\
		$\delta_*$ & [$^\circ$] & $9.5163$ & $3.9528$ & $6.6169$ & $3.5786$ & $3.9910$ & $4.9137$ & $2.7689$ & $4.7237$ \\
		$\alpha_*$ & [$^\circ$] & $-1.3188$ & $3.9297$ & $0.8734$ & $-0.8547$ & $2.1345$ & $0.3825$ & $-0.2585$ & $-0.0014$ \\
		$C_{\text{L,*}}$ & [$-$] & $0.0600$ & $0.0800$ & $0.0400$ & $0.0200$ & $0.0800$ & $0.0400$ & $0.0200$ & $0.0200$ \\
		$C_{\text{D,*}}$ & [$-$] & $0.0006$ & $0.0013$ & $0.0004$ & $0.0001$ & $0.0011$ & $0.0002$ & $0.0001$ & $0.0001$ \\
		$K_{\text{L}}$ & [$\si{kg/m}$] & $0.3135$ & $0.9627$ & $0.6014$ & $0.7220$ & $0.6248$ & $0.6282$ & $0.4471$ & $0.3103$ \\
		$K_{\text{D}}$ & [$\si{kg/m}$] & $0.0031$ & $0.0177$ & $0.0047$ & $0.0023$ & $0.0094$ & $0.0036$ & $0.0010$ & $0.0010$ \\
		$K_{\phi}$ & [$\si{kg}$] & $0.0157$ & $0.0722$ & $0.0376$ & $0.0452$ & $0.0382$ & $0.0479$ & $0.0303$ & $0.0224$ \\
		$K_{\theta}$ & [$\si{kg}$] & $0.0125$ & $0.0562$ & $0.0295$ & $0.0211$ & $0.0226$ & $0.0223$ & $0.0170$ & $0.0175$ \\
		$K_{\psi}$ & [$\si{kg}$] & $0.0002$ & $0.0013$ & $0.0003$ & $0.0001$ & $0.0006$ & $0.0003$ & $0.0001$ & $0.0001$ \\
		$\text{T}_{\text{max}}$ & [$\si{N}$] & $25.0000$ & $25.0000$ & $25.0000$ & $25.0000$ & $25.0000$ & $25.0000$ & $25.0000$ & $25.0000$ \\
		$\text{T}_{\text{min}}$ & [$\si{N}$] & $0.0000$ & $0.0000$ & $0.0000$ & $0.0000$ & $0.0000$ & $0.0000$ & $0.0000$ & $0.0000$ \\
		$\delta_{\text{max}}$ & [$\si{^\circ}$] & $15.0000$ & $15.0000$ & $15.0000$ & $15.0000$ & $15.0000$ & $15.0000$ & $15.0000$ & $15.0000$ \\
		$\delta_{\text{min}}$ & [$\si{^\circ}$] & $15.0000$ & $15.0000$ & $15.0000$ & $15.0000$ & $15.0000$ & $15.0000$ & $15.0000$ & $15.0000$ \\
	\end{tabular}
	\caption{Lower bound, upper bound, baseline and optimal design variables of the tail-sitter. Values above the horizontal line can be chosen by the optimizer. Values below it are fully determined by the ones above it.}
	\label{table:variables}
\end{table*}

Under the assumption that the centre of gravity can be positioned at the midpoint of the symmetry chord, the aircraft carries $0.5\si{kg}$ avionic equipment and two propulsion systems each weighing $0.25\si{kg}$ and is constructed of $100\si{kg/m^3}$ dense EPP foam, the mass characteristics can be readily calculated. This parameterization allows us to compute the aircraft's inertia matrix using the Steiner theorem. The use of AVL enables us to trim the aircraft (i.e., $C_{\theta}(\delta_1,\delta_2)=0$) for a predefined lift coefficient. The drag coefficient $C_{\text{D}}(\alpha)$ is defined as the sum of the zero lift drag coefficient $C_{\text{D},0}$ and the induced drag coefficient, being equal to $k\cdot C_{\text{L}}(\alpha)^2$, with $k^{-1}=\pi\cdot AR\cdot e$ in which $AR$ is the aspect ratio and $e$ the span efficiency factor. $k$ is obtained by means of a polynomial fit using the high fidelity model. 

The operating conditions $\{\alpha_*, \delta_*, C_{\text{L},*}, C_{\text{D}_*}\}$ for which the $\phi$-coefficients ($K_{\text{L}}$, $K_{\text{D}}$, $K_{\phi}$, $K_{\theta}$, and $K_{\psi}$) can determined are chosen such that $(C_{\text{L}}/C_{\text{d}})$ is maximal in trimmed conditions. It is important to note that these coefficients remain constant throughout the trajectory optimization process. The baseline geometry is set to match the centre of the design space and is designed to mimic the commercially available UX5 \cite{wauters:2020b}.

\subsection{Outer loop conceptual design}

Efficiency is the objective that drives the conceptual design, which translates itself to the minimization of the expected value and variance of the time integral over the sum of the thrust required to perform the four manoeuvrers obtained by means of the low-fidelity trajectory optimization, when flown closed loop with optimally tuned control parameters in the high-fidelity environment (Eq.~\ref{eq:CDP}).
\begin{equation}
	\label{eq:l_CDP}
	l_{\text{CDP}}=\text{T}_1+\text{T}_2
\end{equation}

The uncertainty propagation routine employed is Monte Carlo sampling inspired: to obtain the expected value and variance of the time integral, the different missions are flown closed loop in the high fidelity environment a number of times and the metrics are approximated by means of the mean and sample variance respectively. 
\begin{subequations}\label{eq:UP}
	\begin{align}
		L(\vectorstyle{d})&\triangleq \int_0^Tl_{\text{CDP}}(t,\vectorstyle{d}|\vectorstyle{\epsilon})\text{d}t \\
		\mathbb{E}_{\vectorstyle{\epsilon}}\left[L(\vectorstyle{d})\right]& \approx\frac{1}{N}\sum_{i=1}^N L(\vectorstyle{d}) \\
		\mathbb{V}_{\vectorstyle{\epsilon}}\left[L(\vectorstyle{d})\right] &\approx\frac{1}{N}\sum_{i=1}^N \left(L(\vectorstyle{d})-\mathbb{E}_{\vectorstyle{\epsilon}}[L(\vectorstyle{d})]\right)^2
	\end{align}
\end{subequations}

The computational cost that comes with the nested formulation motivates the inclusion of multi-objective \textsl{Bayesian optimization} (BO) in the outer loop \cite{wauters:2022c}. In case the trajectory optimizer is unable to find a feasible path or the uncertainty propagation routine fails to quantify the tracked variance, which might indicate the inability of the design to fly in physical conditions using the reference and tracking strategy provided, the design is denoted as `failed' and treated using the Bayesian classifier. The use of BO within the context of multi-objective design with failed evaluations has been tackled in \cite{wauters:2024}. The framework (without the analytic UP presented in the reference) is used here and briefly detailed in Appendix \ref{sec:BO}. Furthermore, since the UP approach is sampling based, mean and variance will deviate from the exact value, therefore a regressive Gaussian process interpolator is employed, which is also detailed in the Appendix.

To initialize the BO routine, a \textsl{design of experiments} (DoE) is employed by means of a Sobol sequence containing $d({\vectorstyle{d}})\times11-5$ designs. Additionally, 200 iterations can be performed by the optimizer or the routine is halted early if the improvement that can be made (Eq.~\ref{eq:EI}) drops below 0.1\% (scaled by the improvement obtained in the first iteration). A genetic algorithm is employed to optimize the acquisition function.

\subsection{Inner loop trajectory optimization}\label{sec:traj}

The objective that drives the trajectory optimization routine corresponds to a weighted sum of the efficiency motivated thrust needed to execute the manoeuvrers and the snap, which is the fourth derivative of position, and yaw acceleration. The latter two terms are motivated by the differential flatness formulation: in order to minimize the magnitude of the control inputs, it is more efficient to minimize the flat states that most strongly contribute to the signals than constraining them directly. Secondly, the minimization also ensures the smoothness of the trajectories that therefore stabilizes the optimizer and adds to the trackability in the outer loop. However, since minimizing snap and yaw acceleration doesn't permit weighting the relative contribution of thrust and elevon deflection, the thrust is added as an independent term.
\begin{equation}
	\label{eq:l_OCP}
	l_{\text{OCP}}=w_1||\ddddot{\vectorstyle{p}}^{\text{global}}||^2 + w_2||\ddot{\gamma}||^2+w_3(\text{T}_1+\text{T}_2)
\end{equation}

\noindent with $w_1=10^{-8}$, $w_2=10^{-6}$ and $w_2=10^{-5}$. The optimization is performed using SQP as implemented in \MATLAB. To initialize the optimizer, a solution is provided that meets the constraints. This is realized by means of a third order curve for each of the states, which is subsequently translated to the B-spline coefficients, namely a 52 dimensional vector corresponding to four times 13 basis functions from sixth order B-splines with eight knots equally spaced across the time domain.

\subsection{Inner loop control tuning}

In order to ensure that design will be able to perform the manoeuvrers in a physical environment, the gain matrices of the controller $\{\vectorstyle{k}_j^i|j=p,q;i=0,1,2\}$ are adapted with the objective of minimizing the time integral of the objective $l_{\text{FCD}}$. This objective contains on the one hand the ability of the controller to follow the reference trajectory $\vectorstyle{\xi}_*$ and on the other hand the gradient of the control signal. The motivation of the former is self-evident, while the latter ensures that the gain matrices do not take on unrealistic proportions and turn the controller in bang-bang controller.
\begin{equation}
	\label{eq:l_FCD}
	l_{\text{FCP}}=w_4||\vectorstyle{\xi}_*-\vectorstyle{\xi}||^2 + w_5||\dot{\vectorstyle{\upsilon}}||^2
\end{equation}

\noindent with $w_4=0.1$ and $w_5=1.0$. Note that an increase of $w_5/w_4$ leads to smaller gains.

The same uncertainty propagation method is used as has been introduced for the outer loop conceptual design (Eq.~\ref{eq:UP}). Additionally, the same optimization strategy (with the exception of a single-objective instead of a multi-objective treatment) can be employed. 

\begin{landscape}
\begin{table}[t]
	\centering
	\renewcommand{\arraystretch}{1.4}
	\begin{tabular}{r l | r r r r | r r r r | r r r r | r r r r}
		\multicolumn{2}{c}{state} & \multicolumn{4}{c}{cruise}  & \multicolumn{4}{c}{turn}  & \multicolumn{4}{c}{take-off}  & \multicolumn{4}{c}{landing} \\
		& &  BL & D1 & D2 & D3 & BL & D1 & D2 & D3 & BL & D1 & D2 & D3 & BL & D1 & D2 & D3 \\
		\hline
		$\vectorstyle{p}^{(1)}_{\text{in}}$ & [$\si{m}$] & $0.00$ & $0.00$ & $0.00$ & $0.00$ & $0.00$ & $0.00$ & $0.00$ & $0.00$ & $0.00$ & $0.00$ & $0.00$ & $0.00$ & $0.00$ & $0.00$ & $0.00$ & $0.00$ \\
		$\vectorstyle{p}^{(2)}_{\text{in}}$ & [$\si{m}$] & $0.00$ & $0.00$ & $0.00$ & $0.00$ & $-30.64$ & $-40.91$ & $-21.53$ & $-32.41$ & $0.00$ & $0.00$ & $0.00$ & $0.00$ & $0.00$ & $0.00$ & $0.00$ & $0.00$ \\
		$\vectorstyle{p}^{(3)}_{\text{in}}$ & [$\si{m}$] & $100.00$ & $100.00$ & $100.00$ & $100.00$ & $100.00$ & $100.00$ & $100.00$ & $100.00$ & $0.00$ & $0.00$ & $0.00$ & $0.00$ & $0.00$ & $0.00$ & $0.00$ & $0.00$ \\
		$\vectorstyle{q}^{(1)}_{\text{in}}$ & [$\si{rad}$] & $0.00$ & $0.00$ & $0.00$ & $0.00$ & $0.00$ & $0.00$ & $0.00$ & $0.00$ & $0.00$ & $0.00$ & $0.00$ & $0.00$ & $0.00$ & $0.00$ & $0.00$ & $0.00$ \\
		$\vectorstyle{q}^{(2)}_{\text{in}}$ & [$\si{rad}$] & $0.00$ & $0.00$ & $0.00$ & $0.00$ & $-0.95$ & $-0.95$ & $-0.95$ & $-0.95$ & $-3.14$ & $-3.14$ & $-3.14$ & $-3.14$ & $0.00$ & $0.00$ & $0.00$ & $0.00$ \\
		$\vectorstyle{q}^{(3)}_{\text{in}}$ & [$\si{rad}$] & $-0.05$ & $-0.03$ & $-0.07$ & $-0.04$ & $-0.15$ & $-0.09$ & $-0.22$ & $-0.13$ & $0.73$ & $0.81$ & $0.67$ & $0.78$ & $-0.07$ & $-0.05$ & $-0.06$ & $-0.06$ \\
		$\vectorstyle{v}^{(1)}_{\text{in}}$ & [$\si{m/s}$] & $27.50$ & $31.54$ & $23.29$ & $27.95$ & $20.60$ & $23.82$ & $17.23$ & $21.17$ & $5.00$ & $5.00$ & $5.00$ & $5.00$ & $27.00$ & $27.00$ & $27.00$ & $27.00$ \\
		$\vectorstyle{v}^{(2)}_{\text{in}}$ & [$\si{m/s}$] & $0.00$ & $0.00$ & $0.00$ & $0.00$ & $0.00$ & $0.00$ & $0.00$ & $0.00$ & $0.00$ & $0.00$ & $0.00$ & $0.00$ & $0.00$ & $0.00$ & $0.00$ & $0.00$ \\
		$\vectorstyle{v}^{(3)}_{\text{in}}$ & [$\si{m/s}$] & $0.00$ & $0.00$ & $0.00$ & $0.00$ & $0.00$ & $0.00$ & $0.00$ & $0.00$ & $9.99$ & $9.99$ & $10.00$ & $9.99$ & $0.00$ & $0.00$ & $0.00$ & $0.00$ \\
		$\vectorstyle{\Omega}^{(1)}_{\text{in}}$ & [$\si{rad/s}$] & $0.00$ & $0.00$ & $0.00$ & $0.00$ & $0.00$ & $0.00$ & $0.00$ & $0.00$ & $0.00$ & $0.00$ & $0.00$ & $0.00$ & $0.00$ & $0.00$ & $0.00$ & $0.00$ \\
		$\vectorstyle{\Omega}^{(2)}_{\text{in}}$ & [$\si{rad/s}$] & $0.00$ & $0.00$ & $0.00$ & $0.00$ & $0.00$ & $0.00$ & $0.00$ & $0.00$ & $0.00$ & $0.00$ & $0.00$ & $0.00$ & $0.00$ & $0.00$ & $0.00$ & $0.00$ \\
		$\vectorstyle{\Omega}^{(3)}_{\text{in}}$ & [$\si{rad/s}$] & $0.00$ & $0.00$ & $0.00$ & $0.00$ & $0.00$ & $0.00$ & $0.00$ & $0.00$ & $1.30$ & $1.31$ & $1.36$ & $1.32$ & $-0.13$ & $-0.10$ & $-0.09$ & $-0.14$ \\
		\hline
		$\vectorstyle{p}^{(1)}_{\text{out}}$ & [$\si{m}$] & $100.00$ & $100.00$ & $100.00$ & $100.00$ & $0.00$ & $0.00$ & $0.00$ & $0.00$ & $32.67$ & $35.18$ & $29.12$ & $33.81$ & $80.00$ & $96.48$ & $80.00$ & $80.00$ \\
		$\vectorstyle{p}^{(2)}_{\text{out}}$ & [$\si{m}$] & $0.00$ & $0.00$ & $0.00$ & $0.00$ & $30.64$ & $40.91$ & $21.53$ & $32.41$ & $0.00$ & $0.00$ & $0.00$ & $0.00$ & $0.00$ & $0.00$ & $0.00$ & $0.00$ \\
		$\vectorstyle{p}^{(3)}_{\text{out}}$ & [$\si{m}$] & $100.00$ & $100.00$ & $100.00$ & $100.00$ & $100.00$ & $100.00$ & $100.00$ & $100.00$ & $0.00$ & $0.00$ & $0.10$ & $0.00$ & $18.31$ & $20.00$ & $14.67$ & $19.65$ \\
		$\vectorstyle{q}^{(1)}_{\text{out}}$ & [$\si{rad}$] & $0.00$ & $0.00$ & $0.00$ & $0.00$ & $3.14$ & $3.14$ & $3.14$ & $3.14$ & $0.00$ & $0.00$ & $0.00$ & $0.00$ & $0.00$ & $0.00$ & $0.00$ & $0.00$ \\
		$\vectorstyle{q}^{(2)}_{\text{out}}$ & [$\si{rad}$] & $0.00$ & $0.00$ & $0.00$ & $0.00$ & $-0.95$ & $-0.95$ & $-0.95$ & $-0.95$ & $0.00$ & $0.00$ & $0.00$ & $0.00$ & $0.00$ & $0.00$ & $0.00$ & $0.00$ \\
		$\vectorstyle{q}^{(3)}_{\text{out}}$ & [$\si{rad}$] & $-0.05$ & $-0.03$ & $-0.07$ & $-0.04$ & $-0.15$ & $-0.09$ & $-0.21$ & $-0.13$ & $-0.09$ & $-0.08$ & $-0.09$ & $-0.08$ & $-1.54$ & $-1.51$ & $-1.53$ & $-1.46$ \\
		$\vectorstyle{v}^{(1)}_{\text{out}}$ & [$\si{m/s}$] & $27.50$ & $31.54$ & $23.29$ & $27.95$ & $-20.60$ & $-23.82$ & $-17.23$ & $-21.17$ & $27.00$ & $27.00$ & $27.00$ & $27.00$ & $5.00$ & $4.24$ & $5.00$ & $5.00$ \\
		$\vectorstyle{v}^{(2))}_{\text{out}}$ & [$\si{m/s}$] & $0.00$ & $0.00$ & $0.00$ & $0.00$ & $0.00$ & $0.00$ & $0.00$ & $0.00$ & $0.00$ & $0.00$ & $0.00$ & $0.00$ & $0.00$ & $0.00$ & $0.00$ & $0.00$ \\
		$\vectorstyle{v}^{(3)}_{\text{out}}$ & [$\si{m/s}$] & $0.00$ & $-0.00$ & $0.00$ & $0.00$ & $0.00$ & $0.00$ & $0.00$ & $0.00$ & $0.00$ & $0.00$ & $0.00$ & $0.00$ & $5.00$ & $4.24$ & $5.00$ & $5.00$ \\
		$\vectorstyle{\Omega}^{(1)}_{\text{out}}$ & [$\si{rad/s}$] & $0.00$ & $0.00$ & $0.00$ & $0.00$ & $0.00$ & $0.00$ & $0.00$ & $0.00$ & $0.00$ & $0.00$ & $0.00$ & $0.00$ & $0.00$ & $0.00$ & $0.00$ & $0.00$ \\
		$\vectorstyle{\Omega}^{(2)}_{\text{out}}$ & [$\si{rad/s}$] & $0.00$ & $0.00$ & $0.00$ & $0.00$ & $0.00$ & $0.00$ & $0.00$ & $0.00$ & $0.00$ & $0.00$ & $0.00$ & $0.00$ & $0.00$ & $0.00$ & $0.00$ & $0.00$ \\
		$\vectorstyle{\Omega}^{(3)}_{\text{out}}$ & [$\si{rad/s}$] & $0.00$ & $0.00$ & $0.00$ & $0.00$ & $0.00$ & $0.00$ & $0.00$ & $0.00$ & $-0.29$ & $-0.29$ & $-0.31$ & $-0.29$ & $-1.43$ & $-1.20$ & $-1.26$ & $-1.43$ \\
		\hline
		T & [$\si{s}$] & $3.63$ & $3.17$ & $4.29$ & $3.57$ & $4.67$ & $5.39$ & $3.92$ & $4.81$ & $2.17$ & $2.29$ & $1.95$ & $2.21$ & $4.33$ & $5.37$ & $4.35$ & $4.36$ \\
		$l_{\text{OCP}}$ & [$\si{N}$] & $17.25$ & $9.65$ & $29.26$ & $13.29$ & $837.49$ & $533.75$ & $1201.10$ & $648.97$ & $2394.40$ & $2167.00$ & $2892.50$ & $2103.50$ & $284.77$ & $219.69$ & $411.64$ & $201.96$ \\
	\end{tabular}
	\caption{Lower bound, upper bound, baseline and optimal design variables of the tail-sitter. Values above the horizontal line can be chosen by the optimizer. Values below it are fully determined by the ones above it.}
	\label{table:missions}
\end{table}
\end{landscape}

\begin{figure}[b]
	\centering
	\vspace{-5mm}
	\includegraphics[trim=0 0 0 0,clip,width=0.99\linewidth]{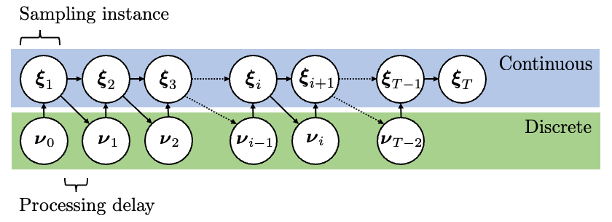}
	\caption{Reality emulating feedback control updating approach.}	
	\label{fig:control_timing}
\end{figure}

Additionally, to further mimic real-life application, the simulation of the closed-loop flight is performed by means of what can be described as a first-order hold delayed feedback control updating approach (Fig.~\ref{fig:control_timing}). This is the behaviour that can be expected due to the delay in sensing and computing of the online unit. In this work we employ a $100\si{Hz}$ sampling rate.

\subsection{Extended Design System Matrix}

When the different components of the tail-sitter dynamics discussed in sec.~\ref{sec:appl}, along with the problem implementation discussed in this section are plugged in the framework introduced in sec.~\ref{sec:codesign}, the full first-time-right dynamic design approach of a tail-sitter is obtained. This is presented in Fig.~\ref{fig:XDSM2}. The acronyms used are summarized in the nomenclature at the end of the document.

\subsection{Results and discussion}

\begin{figure}[t]
	\centering
	\includegraphics[trim=0 25 0 25,clip,width=0.95\linewidth]{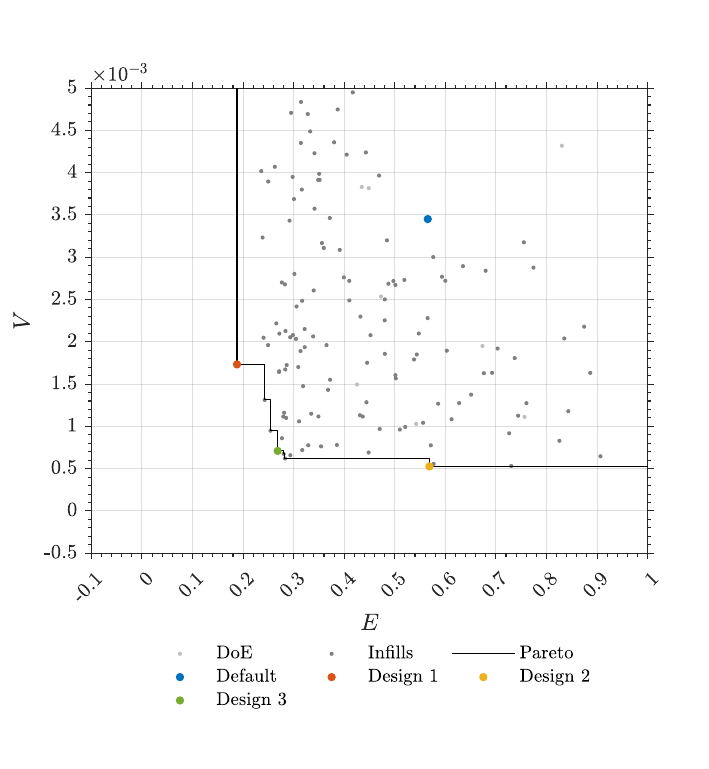}
	\caption{Pareto front of Tail-sitter designs according to the robust dynamic design framework}	
	\label{fig:pareto}
	\vspace{-5mm}
\end{figure}

The result of the first-time-right dynamic design architecture is a Pareto front of tail-sitter designs, trajectories and control settings that outperform the baseline design both in mean and variance of the energy expenditure required to perform the four manoeuvres (Fig.~\ref{fig:pareto}). The corresponding design variables and controller settings are presented in Table~\ref{table:variables}. Furthermore, the boundary conditions of the different missions are presented in Table~\ref{table:missions}. 

\begin{figure*}[h!]
	\vspace{-7mm}
	\centering
	\begin{subfigure}{0.24\linewidth}
		\centering
		\includegraphics[trim=0 0 0 0,clip,width=0.99\linewidth]{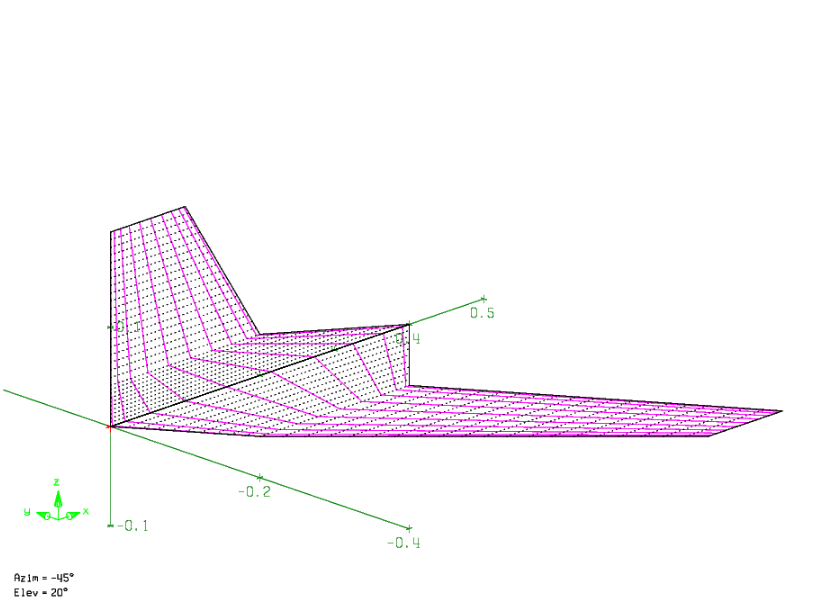}
		\caption{Lower bound design}
	\end{subfigure}
	\begin{subfigure}{0.24\linewidth}
		\centering
		\includegraphics[trim=0 0 0 0,clip,width=0.99\linewidth]{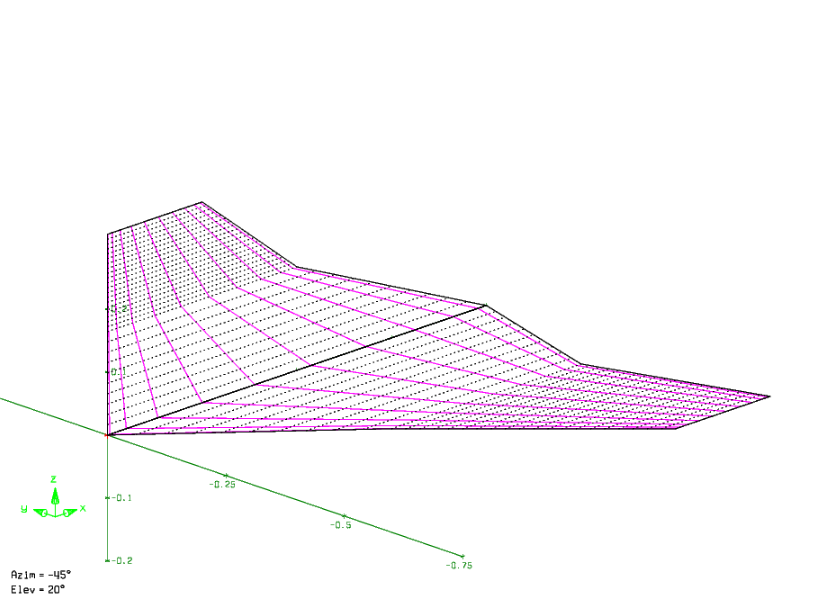}
		\caption{Upper bound design}
	\end{subfigure}
	\begin{subfigure}{0.24\linewidth}
		\centering
		\includegraphics[trim=0 0 0 0,clip,width=0.99\linewidth]{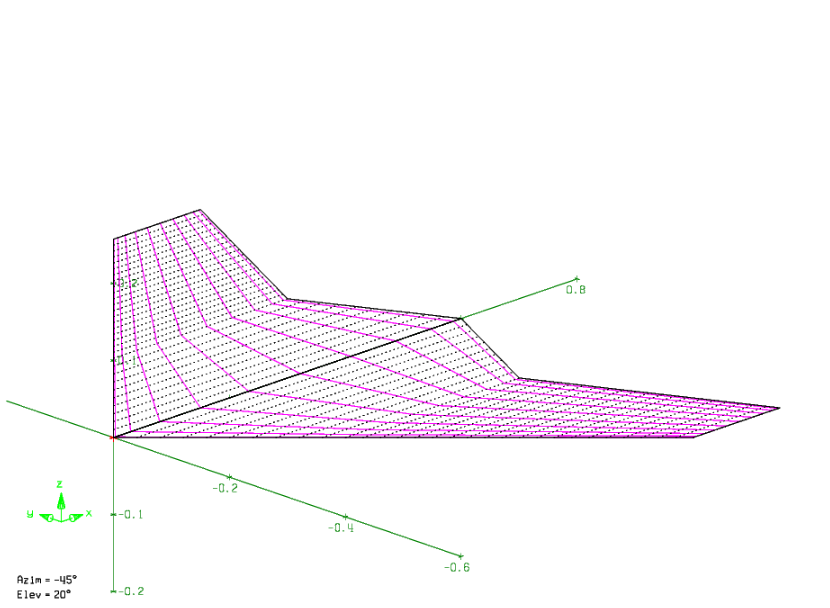}
		\caption{Baseline design}
	\end{subfigure}
	\begin{subfigure}{0.24\linewidth}
		\centering
		\includegraphics[trim=0 0 0 0,clip,width=0.99\linewidth]{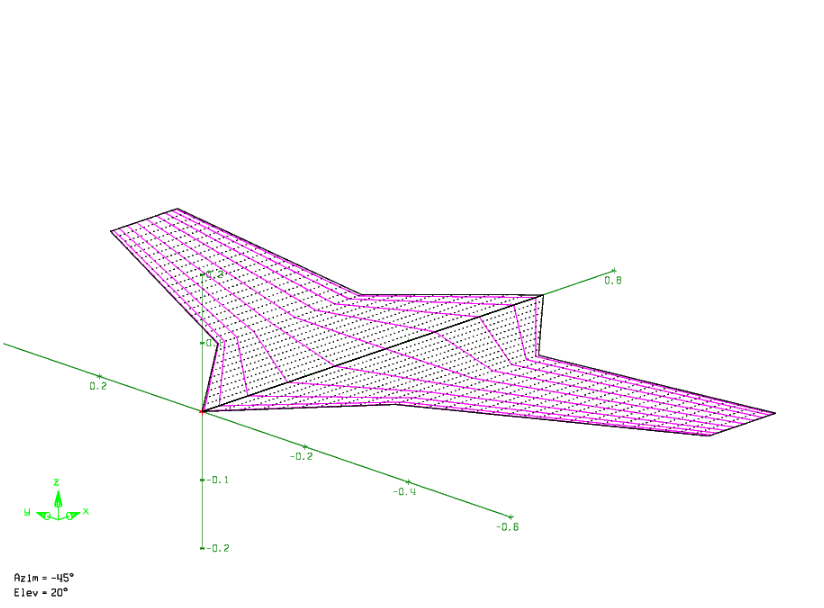}
		\caption{Design 1}
	\end{subfigure}
	\begin{subfigure}{0.24\linewidth}
		\centering
		\includegraphics[trim=0 0 0 0,clip,width=0.99\linewidth]{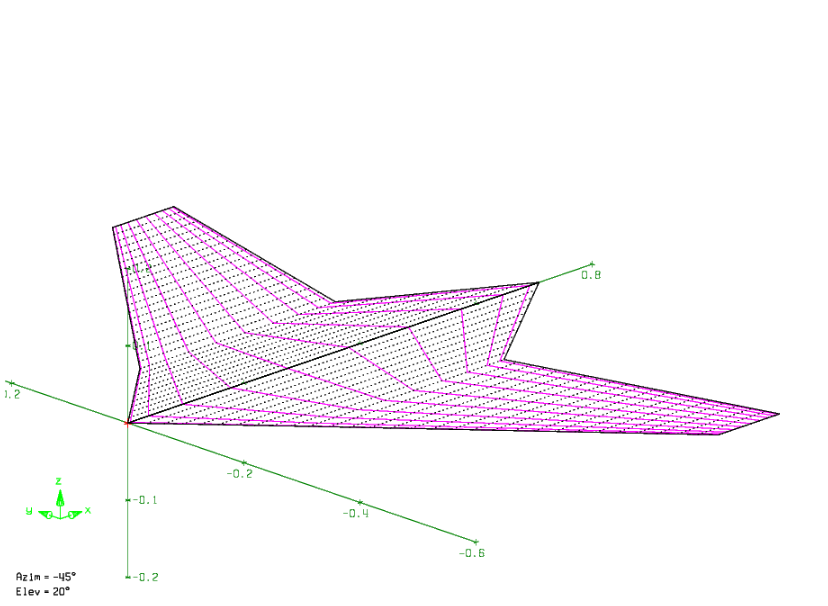}
		\caption{Design 2}
	\end{subfigure}
	\begin{subfigure}{0.24\linewidth}
		\centering
		\includegraphics[trim=0 0 0 0,clip,width=0.99\linewidth]{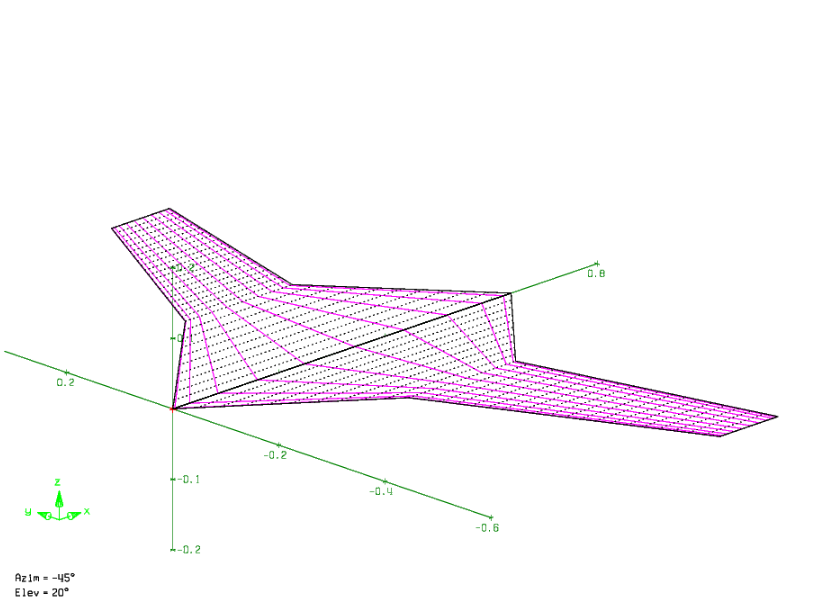}
		\caption{Design 3}
	\end{subfigure}
	\begin{subfigure}{0.24\linewidth}
		\centering
		\includegraphics[trim=0 0 0 0,clip,width=0.99\linewidth]{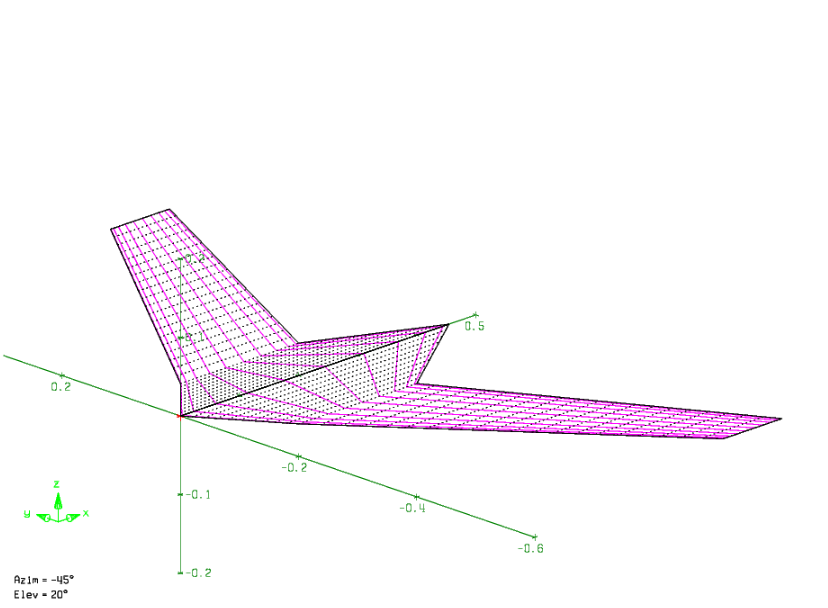}
		\caption{Ablation design 1}
	\end{subfigure}
	\begin{subfigure}{0.24\linewidth}
		\centering
		\includegraphics[trim=0 0 0 0,clip,width=0.99\linewidth]{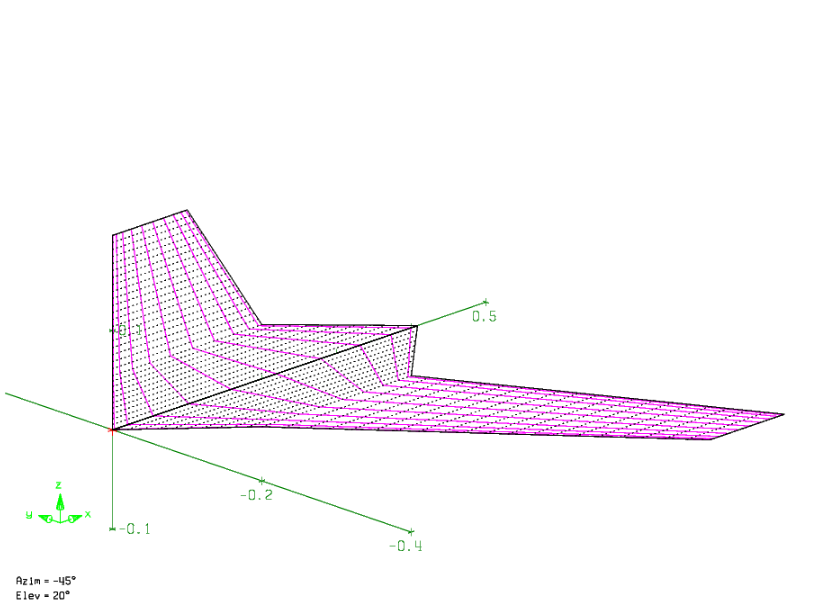}
		\caption{Ablation design 2}
	\end{subfigure}
	\caption{Planform visualization of the designs presented in Table~\ref{table:variables}.}	 
	\label{fig:designs}
\end{figure*}

In comparison with the default design, it can be observed that the designs are more slender (higher aspect ratio), characterized by smaller chord lengths, larger spans and a smaller sweeps. However, in order to ensure that the control authority is maintained, the control surface sizes are increased (Fig.~\ref{fig:designs}). In regards to the control tuning, it can be noted that the integrator term of the outer loop position design is set to zero, while large values are obtained for the integrator term of the inner loop attitude controller. This indicates that the reality gap is completely accounted for by the inner controller.

\begin{figure*}[t]
	\centering
	\includegraphics[trim=40 15 40 15,clip,width=0.99\linewidth]{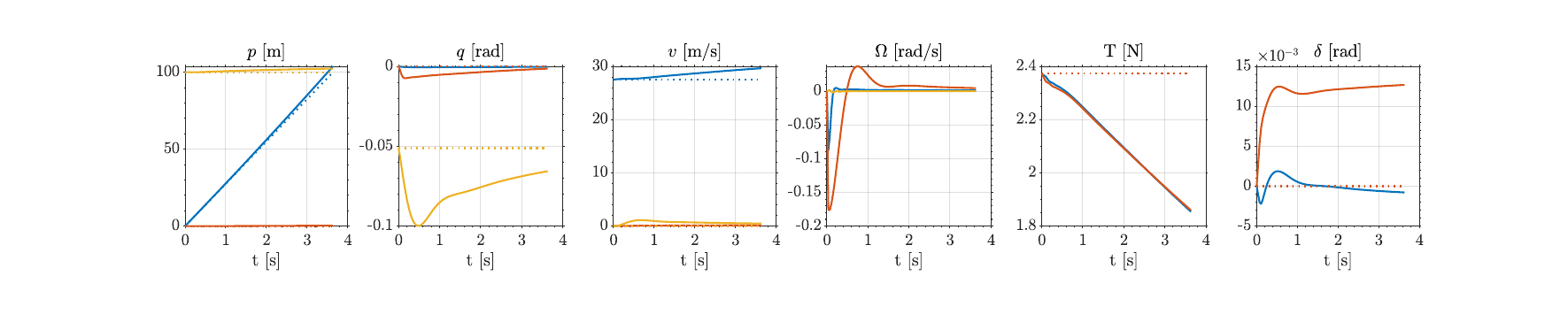}
	\includegraphics[trim=40 15 40 15,clip,width=0.99\linewidth]{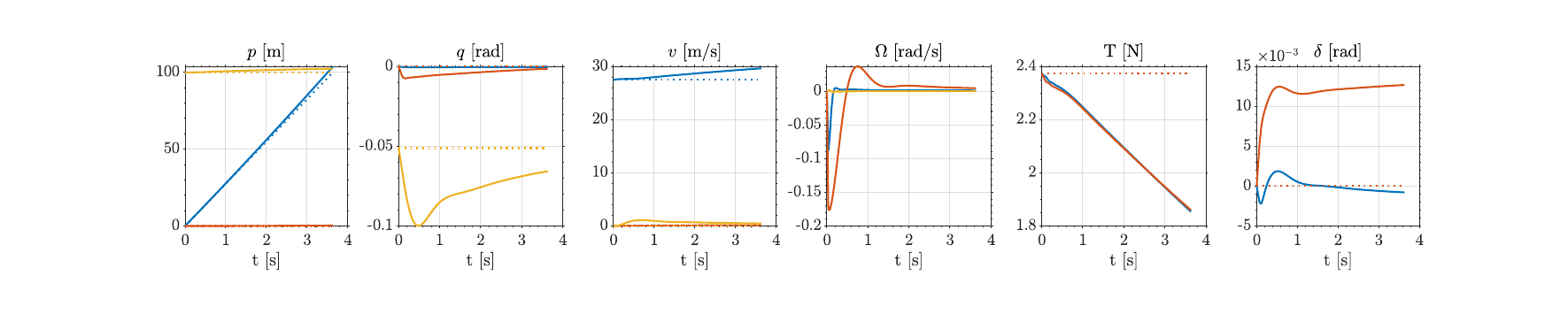}
	\includegraphics[trim=40 15 40 15,clip,width=0.99\linewidth]{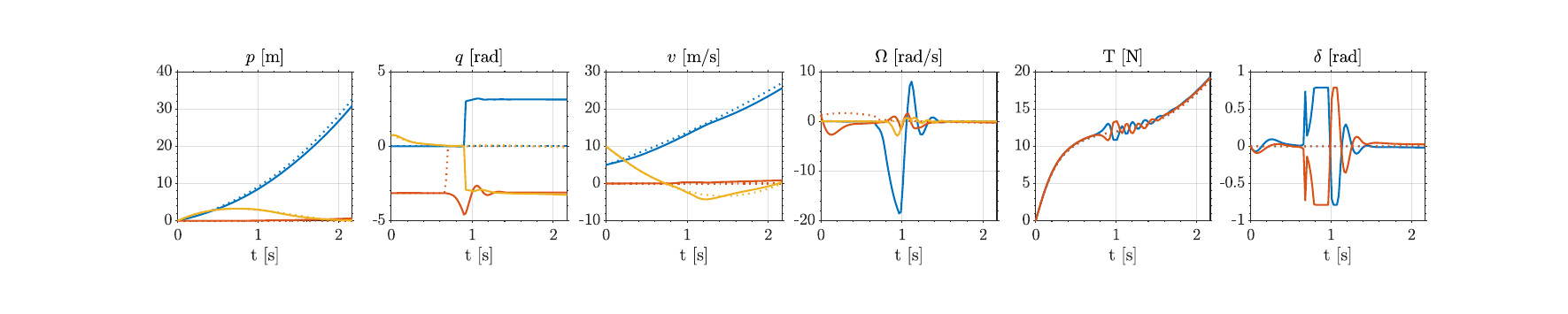}
	\includegraphics[trim=40 15 40 15,clip,width=0.99\linewidth]{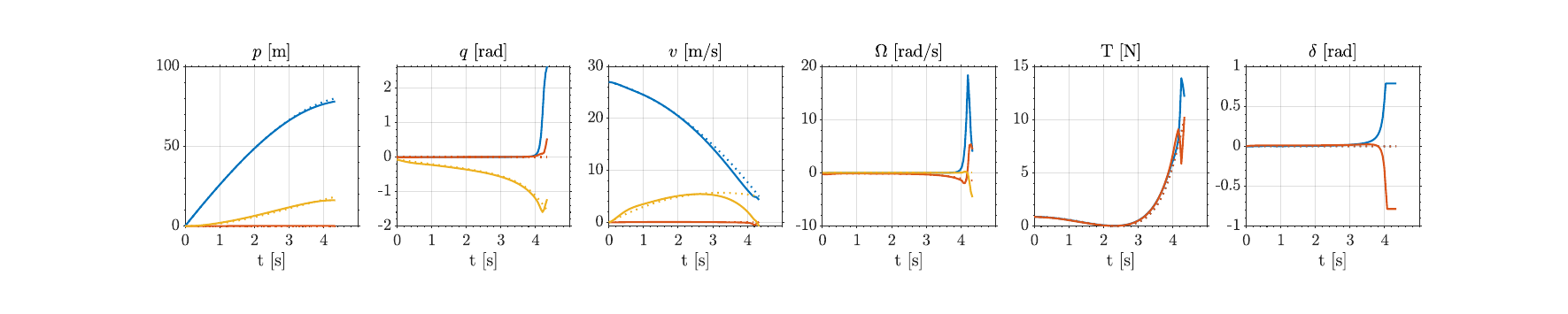}
	\label{fig:resultsBaseline}
	\caption{Comparison of the (low-fidelity) reference trajectory (dashed line) and the (high-fidelity) closed-loop trajectory (full line) of the baseline design for respectively from top to bottom cruise, turn, take-off and landing manoeuvre. First, second and third arguments are respectively displayed in blue, red and yellow.}
\end{figure*}

\begin{figure*}[t]
	\centering
	\includegraphics[trim=40 15 40 15,clip,width=0.99\linewidth]{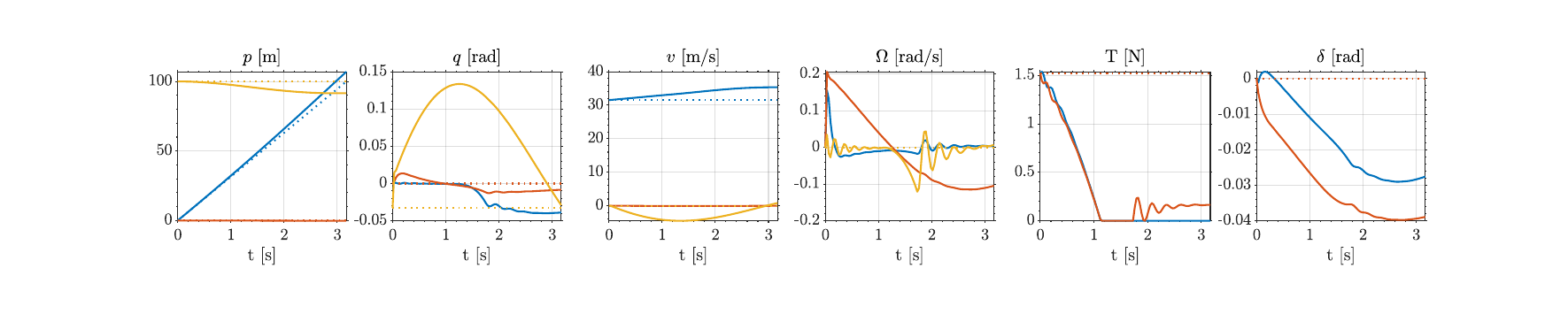}
	\includegraphics[trim=40 15 40 15,clip,width=0.99\linewidth]{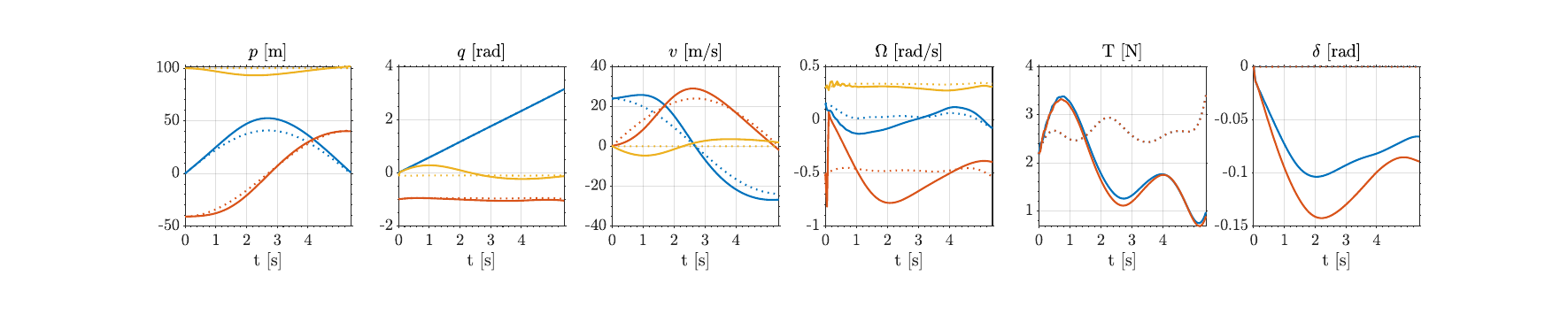}
	\includegraphics[trim=40 15 40 15,clip,width=0.99\linewidth]{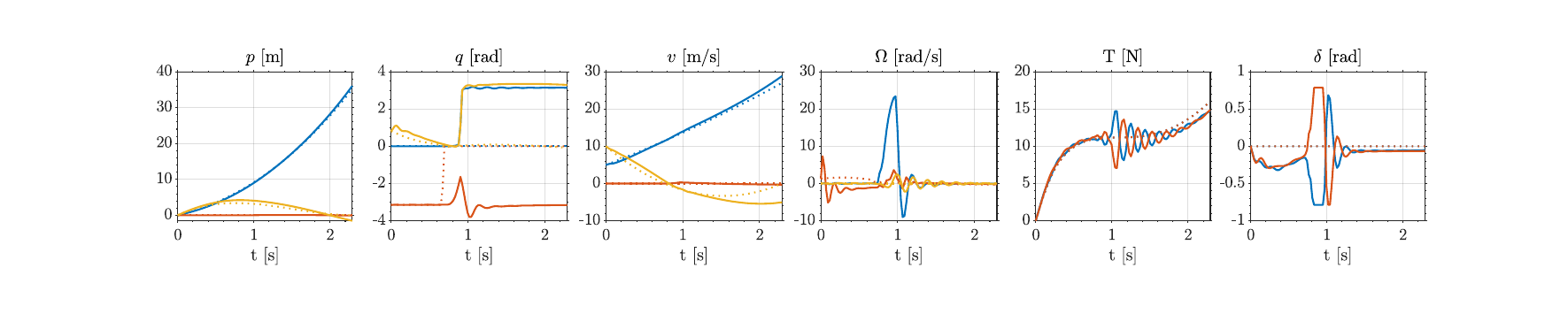}
	\includegraphics[trim=40 15 40 15,clip,width=0.99\linewidth]{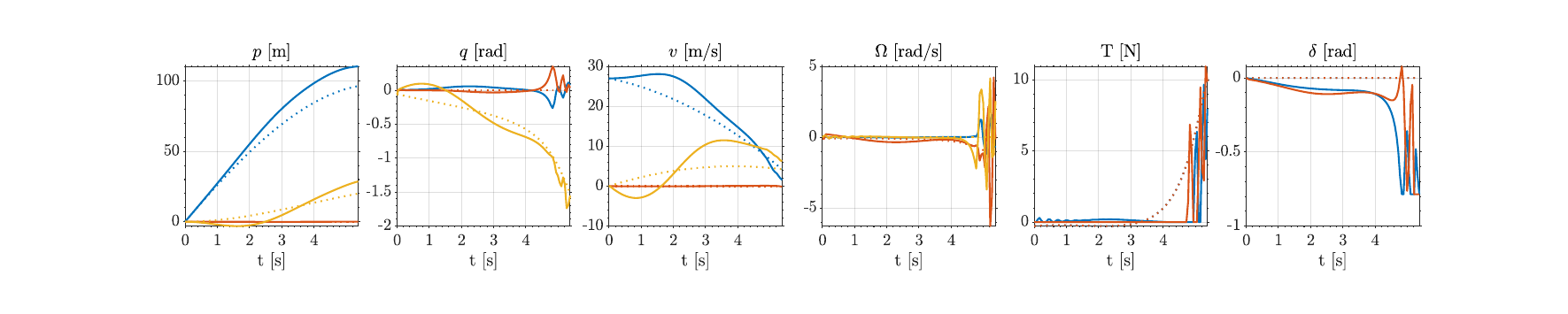}
	\label{fig:resultsDesign1}
	\caption{Comparison of the (low-fidelity) reference trajectory (dashed line) and the (high-fidelity) closed-loop trajectory (full line) of design 1 for respectively from top to bottom cruise, turn, take-off and landing manoeuvre. First, second and third arguments are respectively displayed in blue, red and yellow.}
\end{figure*}

\begin{figure*}[t]
	\centering
	\includegraphics[trim=40 15 40 15,clip,width=0.99\linewidth]{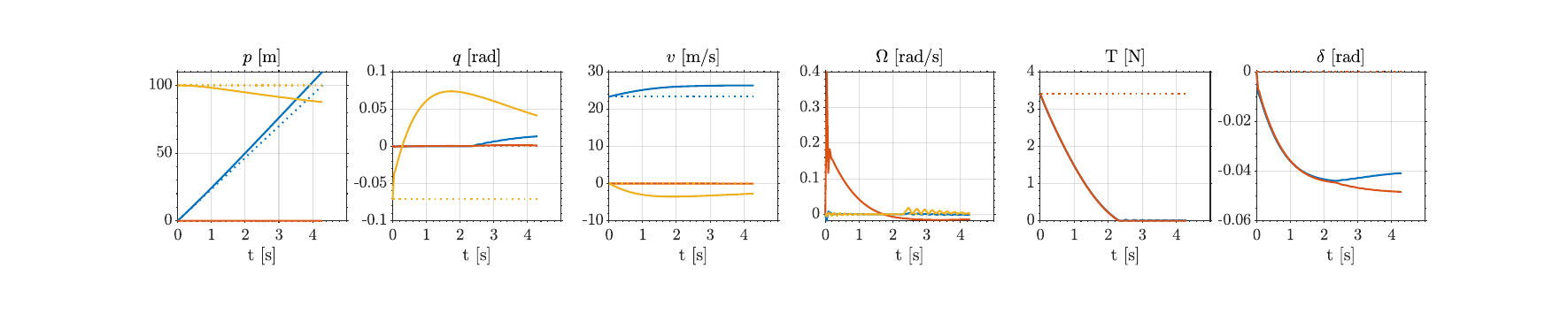}
	\includegraphics[trim=40 15 40 15,clip,width=0.99\linewidth]{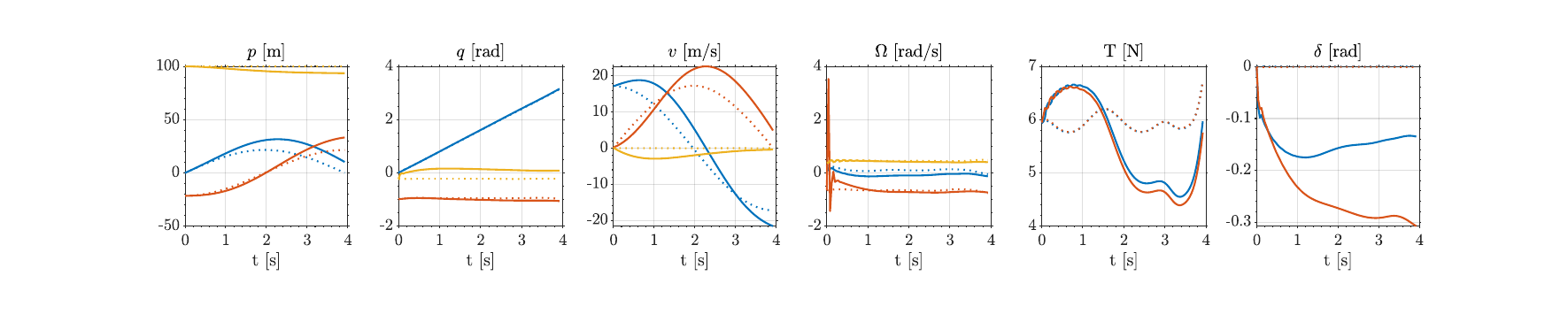}
	\includegraphics[trim=40 15 40 15,clip,width=0.99\linewidth]{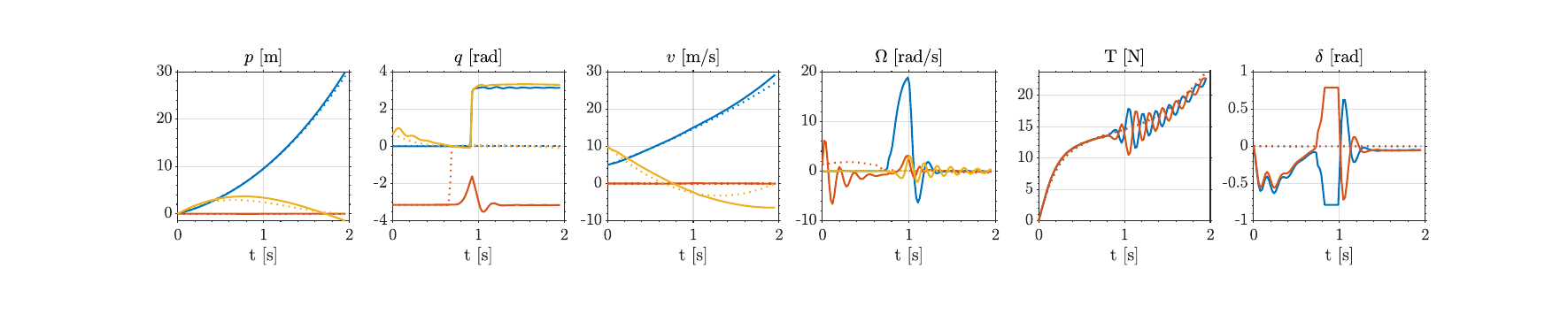}
	\includegraphics[trim=40 15 40 15,clip,width=0.99\linewidth]{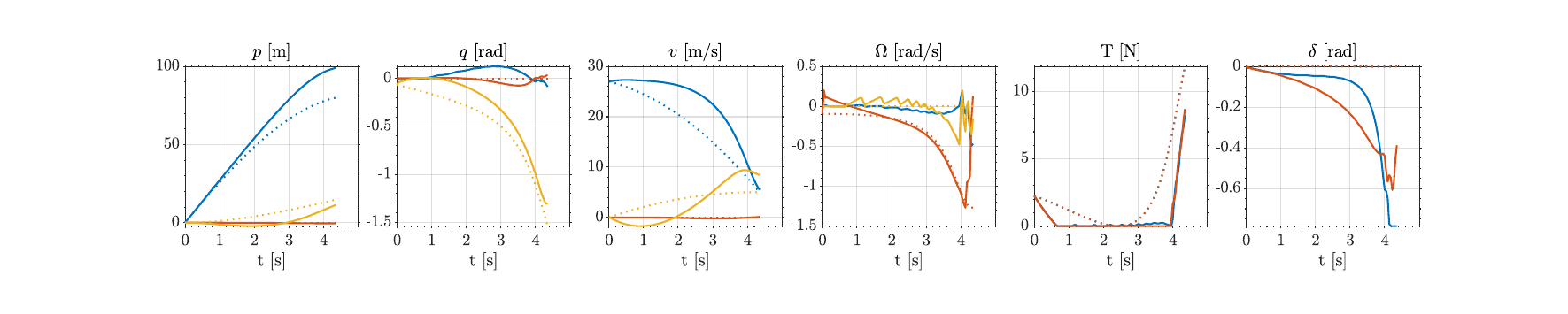}
	\label{fig:resultsDesign2}
	\caption{Comparison of the (low-fidelity) reference trajectory (dashed line) and the (high-fidelity) closed-loop trajectory (full line) of design 2 for respectively from top to bottom cruise, turn, take-off and landing manoeuvre. First, second and third arguments are respectively displayed in blue, red and yellow.}
\end{figure*}

\begin{figure*}[t]
	\centering
	\includegraphics[trim=40 15 40 15,clip,width=0.99\linewidth]{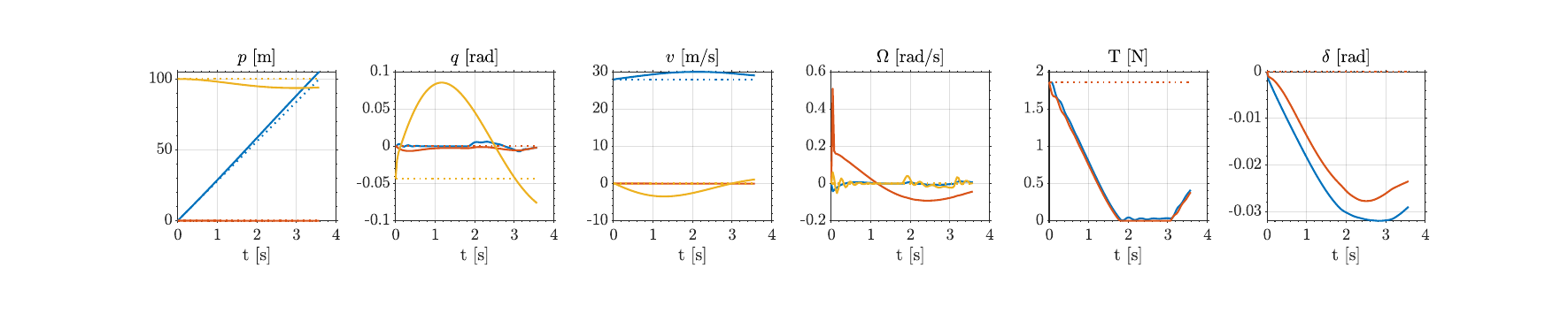}
	\includegraphics[trim=40 15 40 15,clip,width=0.99\linewidth]{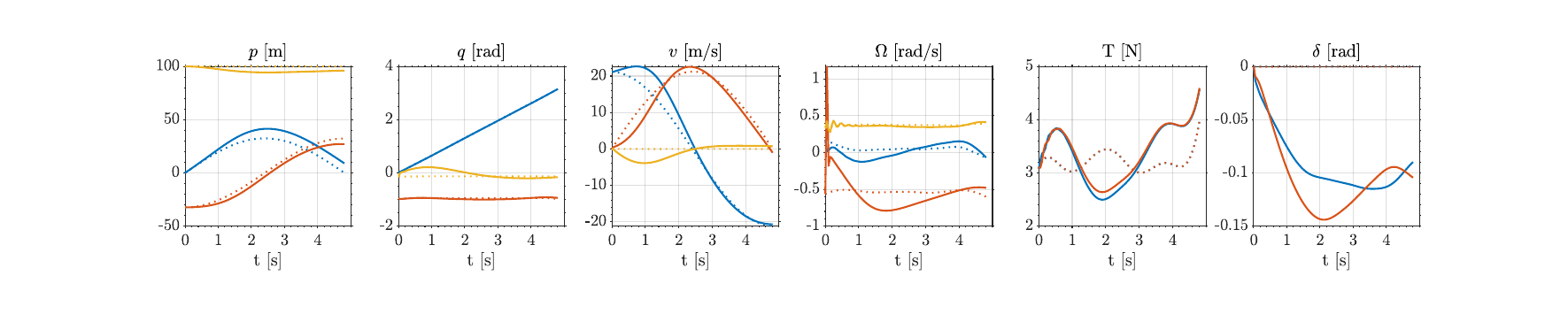}
	\includegraphics[trim=40 15 40 15,clip,width=0.99\linewidth]{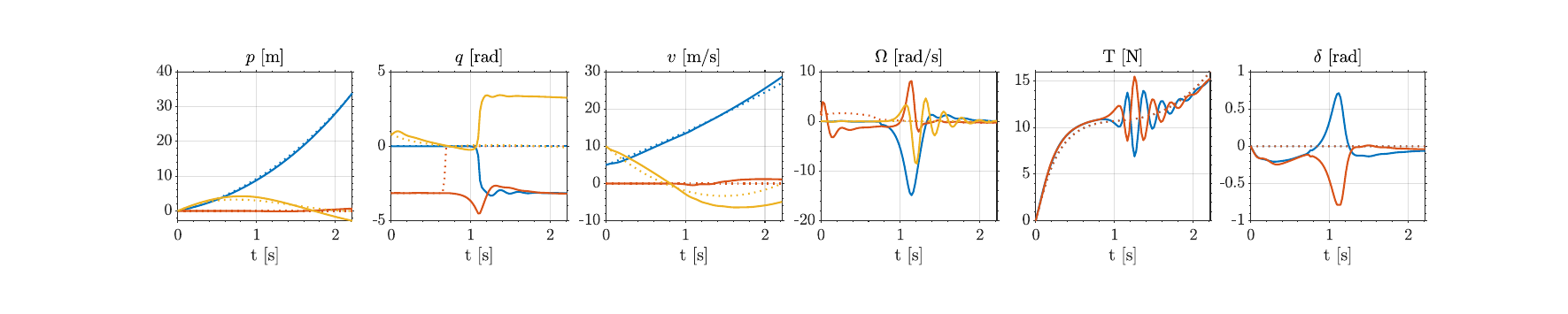}
	\includegraphics[trim=40 15 40 15,clip,width=0.99\linewidth]{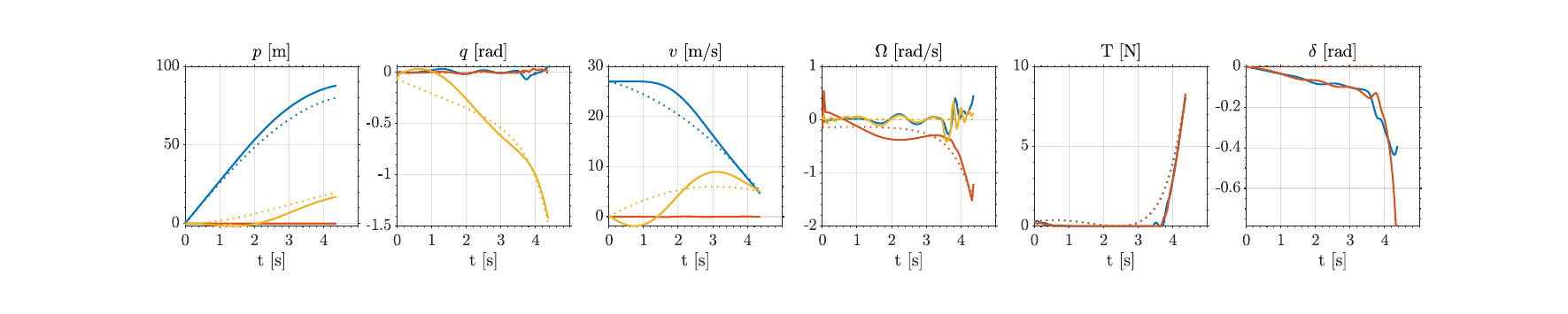}
	\label{fig:resultsDesign3}
	\caption{Comparison of the (low-fidelity) reference trajectory (dashed line) and the (high-fidelity) closed-loop trajectory (full line) of design 3 for respectively from top to bottom cruise, turn, take-off and landing manoeuvre. First, second and third arguments are respectively displayed in blue, red and yellow.}
\end{figure*}

In Fig.~\ref{fig:resultsBaseline},~\ref{fig:resultsDesign1},~\ref{fig:resultsDesign2} and~\ref{fig:resultsDesign3} the four missions flown by respectively the baseline, design 1, design 2 and design 3 are presented by comparing the reference trajectory obtained by the flatness-based trajectory optimization (sec.~\ref{sec:forward} and~\ref{sec:traj}) with the trajectory flown in the reality-emulator using optimally tuned control parameters (sec.~\ref{sec:AVL} and~\ref{sec:tuning}). For these figures, the inherent stochasticity of the reality-emulator is turned off to permit a qualitative discussion. It can be observed that overall the designs are able to fly the missions successfully by examining the position plots on the left. However, when examining the attitude plots, specifically the take-off mission, an abrupt $\pi$ can be observed, which translates itself in abrupt response of $\vectorstyle{\Omega}$, $\vectorstyle{T}$ and $\vectorstyle{\delta}$. Although the controller can account for this unphysical jump, it is nonetheless a numerical hiccup that must be accounted for. This is caused by the flatness transform which includes a $\tan^{-1}$. Restricting the outcome between from $[0,2\pi]$ to $[-\pi,\pi]$ simply moves the problem elsewhere and therefore requires a more subtle treatment. We leave this as future work, as its presence in the current work does not undermine the results presented here. Finally, the landing manoeuvre shows violent behaviour near the end of the mission. This is caused by the reaching the workable limit of the model since it does not include prop wash to ensure the effectiveness of the elevons. Extending the model is currently being examined.

\subsection{Ablation study}

To further illustrate the effectiveness of DAIMYO, two ablation studies are performed that correspond to a nested co-design without the reality-emulator and a more conventional static design approach. In case of the former, this corresponds to the design approach presented in \cite{wauters:2023}. The results are added to Table~\ref{table:variables}. It can be noted that the gain matrices are left blank. This corresponds to failure of the control tuning approach to find suitable gains to fly the missions in the reality-emulator. As such, this would imply the inability of the design to be employed in real-life.

\section{Conclusion and outlook}

In this paper, we present a dynamic design architecture that aims to enhance the multi-disciplinary design of dynamical systems in a practical and reliable manner. Our approach adopts a nested concurrent design formulation that combines both conceptual and control design aspects. The outer loop of our framework employs a fail-safe Bayesian optimization process with constraints to optimize performance metrics based on conceptual design parameters. Simultaneously, the inner loop utilizes a trajectory optimization framework based on differential flatness principles, along with an aerodynamic solver, to determine an energy-efficient flight path.

One significant aspect of our approach is the incorporation of a reality emulator to assess the impact of the reality gap, which represents the mismatch between our simplified model and actual operational conditions. This assessment allows us to tune a feedback controller that ensures the safe operation of the designs once deployed. 

To illustrate the effectiveness of our approach, we applied it to the challenging problem of tail-sitter development, specifically focusing on optimizing Dubins path manoeuvrers. The iterative scheme led to a design that demonstrates both energy efficiency and the potential for strong real-world performance.

While this design is not necessarily the most optimal design according to the model, it is a design for which we can expect that its operation in reality, as defined in our design routine, will perform adequately. While a full validation of the framework would require an experimental implementation, which we leave for future work, the inclusion of the emulator to mimic the operational conditions nonetheless proves the effectiveness of the scheme.

Future work is currently directed in three directions. First, the numerical hiccup occurring during the take-off manoeuvre is examined. Second, the differentally flatness of a tail-sitter model that includes prop wash is being researched. Third, alternative strategies, which on the one hand could nest the feedforward and feedback control design and on the other hand would permit a more simultaneous treatment of the design problem are being developed. 

\section*{Acknowledgment}

The authors gratefully acknowledge the funding by the Research Foundation-Flanders (FWO), Belgium, through the Junior Post-Doctoral fellowship of Jolan Wauters and the Energy Transition Fund (ETF) project BORNE of FPS Economy. The computational resources (Stevin Supercomputer Infrastructure) and services used in this work were provided by the VSC (Flemish Supercomputer Center), funded by Ghent University, Belgium, FWO and the Flemish Government — department EWI, Belgium.

\appendices
\section{Differential flatness of tail-sitter UAVs}\label{sec:flat}

The differential flatness of the low-fidelity tail-sitter model with negligible lateral force was established by \cite{tal:2023}. Here we provide additional details about the model and the associated implications of its flatness. In this appendix we adopt the convention of accentuating vectorial physical quantities as a function of the chosen frame of reference. Note that therefore we can drop the super scripted annotations from section \ref{sec:rbd}. Starting with a vector $\vectorstyle{w}\in\mathbb{R}^3$ in the global frame we define $\vectorstyle{w}' = \mathrm{rot}^\top_z(\psi) \vectorstyle{w}$, $\vectorstyle{w}'' = \mathrm{rot}^\top_x(\phi) \vectorstyle{w}'$ and $\vectorstyle{w}^{\text{local}}=\vectorstyle{w}''' = \mathrm{rot}^\top_y(\theta) \vectorstyle{w}''$ where $\mathrm{rot}_\alpha,\alpha\in\{x,y,z\}$ represents the rotation matrix about the axis $\alpha$.

To establish the differential flatness of the model, first a flat output $\vectorstyle{\sigma}$ needs be determined. Intuitively, the flat output can be interpreted as a minimal dynamical representation of any feasible state-input trajectory. Put differently, the system is flat if there exist differential operators $\Xi$ and $\Upsilon$ so that
\begin{subequations}
	\begin{align}
		\vectorstyle{\xi} &= \Xi[\vectorstyle{\sigma}] = \Xi(\vectorstyle{\sigma},\dot{\vectorstyle{\sigma}},\ddot{\vectorstyle{\sigma}},\dots)\\
		\vectorstyle{\upsilon} &= \Upsilon[\vectorstyle{\sigma}] = \Upsilon(\vectorstyle{\sigma},\dot{\vectorstyle{\sigma}},\ddot{\vectorstyle{\sigma}},\dots)
	\end{align}
\end{subequations}

Tal showed that the tail-sitter is flat with flat output, $\vectorstyle{\sigma}=(x,y,z,\psi)$. In the remainder of this appendix, the differential flatness is demonstrated. Deriving $\vectorstyle{v}$ and $\dot{\vectorstyle{v}}$ as functions of the derivatives of $\vectorstyle{p}$ is a straightforward calculation. Expressing the angle vector poses a challenge. From the equation for $\dot{\vectorstyle{v}}$, we can determine the force vector $\vectorstyle{f}$ in global coordinates. With the knowledge of $\psi$, we can calculate $\vectorstyle{f}'$. Following the geometric aspects of the problem and assuming the absence of a second force component in $\vectorstyle{f}'''$, it becomes evident that $\vectorstyle{f}'$ must lie within the roll plane. In other words we can derive the roll angle $\phi$
\begin{equation}
	\phi = -\arctan \frac{\vectorstyle{f}_y'}{\vectorstyle{f}_z'}
	\label{eq:phi}
\end{equation}

Upon determining the roll angle, we can compute $\vectorstyle{f}''$ and $\vectorstyle{v}''$. Subsequently, we can solve for both $\theta$ and $T$ using the equation (\ref{eq:forces}), where $\text{T}$ represents the total thrust force, which is equal to $\text{T}_1 + \text{T}_2$.
\begin{subequations}
	\begin{align}
		\theta =&\: \arctan \frac{-\vectorstyle{f}_z''-K_{\text{L}} \vectorstyle{v}_z''\|\vectorstyle{v}\|}{\vectorstyle{f}_x'' + K_{\text{L}} \vectorstyle{v}_x''\|\vectorstyle{v}\|} \\
		\text{T} =&\: (-\vectorstyle{f}_z''-K_{\text{D}} \vectorstyle{v}_z''\|\vectorstyle{v}\|) \sin(\theta)  \nonumber\\
		&+ (\vectorstyle{f}_x''+K_{\text{D}} \vectorstyle{v}_x''\|\vectorstyle{v}\|) \cos(\theta)
	\end{align}
	\label{eq:T}
\end{subequations}

By utilizing the angles, we can apply the kinematics of the problem to deduce $\vectorstyle{\Omega}'''$. With knowledge of $\vectorstyle{\Omega}'''$ and utilizing equation (\ref{eq:rbd}), we can compute $\vectorstyle{\tau}'''$. This value is associated with the control surface deflection and provides the means to determine how the total thrust force is distributed across the propulsion system.
\begin{subequations}
	\begin{align}
		\text{T}_1 &= \frac{1}{2}\text{T} + \frac{1}{2l_y^{\text{T}}}\left(\frac{K_{\psi}}{K_{\phi}}\vectorstyle{\tau}_x''' - \vectorstyle{\tau}_z'''\right)\\
		\text{T}_2 &= \frac{1}{2}\text{T} - \frac{1}{2l_y^{\text{T}}}\left(\frac{K_{\psi}}{K_{\phi}}\vectorstyle{\tau}_x''' - \vectorstyle{\tau}_z'''\right)\\
		\delta_1 &= \frac{1}{2\vectorstyle{v}_x'''\|\vectorstyle{v}\|}\left(\frac{\vectorstyle{\tau}_y'''}{K_{\theta}}-\frac{\vectorstyle{\tau}_x'''}{K_{\phi}}\right) \\
		\delta_2 &= \frac{1}{2\vectorstyle{v}_x'''\|\vectorstyle{v}\|}\left(\frac{\vectorstyle{\tau}_x'''}{K_{\phi}}+\frac{\vectorstyle{\tau}'''_y}{K_{\theta}}\right)
	\end{align}
	\label{eq:invDyn}
\end{subequations}

Consequently, we have derived equations expressing $\vectorstyle{\xi}$ and $\vectorstyle{\upsilon}$ in terms of the derivatives of $\vectorstyle{\sigma}$.

\section{Bayesian optimization}\label{sec:BO}

We briefly present the cornerstones of \textsl{Bayesian optimization} (BO), namely \textsl{Gaussian process interpolation} (GPI) and the acquisition function \textsl{expected improvement} (EI), along with a number of extensions that enable the application of the DAIMYO architecture on the tail-sitter design case. A more extensive discussion on the subject is presented in \cite{rasmussen:2006,shahriari:2016}. 

The incentive of BO is finding the global optimizer of a black-box function $\text{y}:\mathcal{X}\mapsto\mathbb{R}$ in a data efficient manner. $\mathcal{X}$ has been introduced to denote the design space, most often a confined subset of $\mathbb{R}^{d({\mathcal{X}})}$, in which $d({\mathcal{X}})$ is used to denote the dimensionality of the input space and, based on the state-of-practice, is taken not larger than 20 \cite{keane:2006}.
\begin{equation}
	\vectorstyle{x}_{*}=\arg\min_{\substack{\vectorstyle{x}\in\mathcal{X}}}\:\text{y}(\vectorstyle{x})
\end{equation}

BO builds on the idea of replacing the expensive objective function by a cheap-to-sample surrogate model, most often a \textsl{Gaussian process interpolator} (GPI), which can be used to determine the next design to be evaluated by means of an acquisition function \cite{shahriari:2016}.

\subsection{Gaussian Process Interpolation}

A Gaussian process $\mathcal{y}$ can be understood as a distribution over functions, such that $\mathcal{y}(\vectorstyle{x})$ follows a joint Gaussian distribution at any set of points $\{\vectorstyle{x}_i|i=1,...,n\}$. Therefore, the stochastic process becomes fully defined by means of a mean function $m(\vectorstyle{x}):\mathcal{X}\mapsto\mathbb{R}$ and a covariance function $\psi(\vectorstyle{x},\vectorstyle{x}'|\vectorstyle{\theta}):\mathcal{X}\times\mathcal{X}\mapsto\mathbb{R}$, which is in turn a function of hyperparameters $\vectorstyle{\theta}$. As is customary in the aerospace community, we employ the Mat\'{e}rn covariance function \cite{palar:2020}, of which a full description is given in Rasmussen \& Williams, chapter 4 \cite{rasmussen:2006}.

In order to obtain the hyperparameters, such that the Gaussian process best reproduces the data, use can be made of a point estimate in the shape of the \textsl{maximum concentrated log marginal likelihood estimation} (MLE):
\begin{equation}\label{eq:MLE}
	\vectorstyle{\theta}_*=\arg\max_{\vectorstyle{\theta}}\:-n(\matrixstyle{S}^{\vectorstyle{i}})\log\sigma^2-\log|\matrixstyle{\Psi}|
\end{equation}
with $n(\matrixstyle{S}^{\vectorstyle{i}})$ the number of elements in the evaluated data set $\matrixstyle{S}^{\vectorstyle{i}}\triangleq\{\matrixstyle{X},\vectorstyle{y}\}$ the evaluated data, $\matrixstyle{\Psi}$ the correlation matrix defined by $\matrixstyle{\Psi}^{(i,j)}=\psi(\vectorstyle{x}_i,\vectorstyle{x}_j)$ and $\sigma^2$ the process variance, obtained by means of a \textsl{Generalized Least Squares} (GLS) approach such that
\begin{subequations}
	\begin{align}
		\vectorstyle{\beta}=&\:(\matrixstyle{F}^\top\matrixstyle{\Psi}^{-1}\matrixstyle{F})^{-1}\matrixstyle{F}^\top\matrixstyle{\Psi}^{-1}\textbf{y} \\ 
		\sigma^2=&\:n(\matrixstyle{S}^{\vectorstyle{i}})^{-1}(\vectorstyle{y}-\matrixstyle{F}\vectorstyle{\beta})^\top\matrixstyle{\Psi}^{-1}(\vectorstyle{y}-\matrixstyle{F}\vectorstyle{\beta})
	\end{align}
\end{subequations}

\noindent with $\matrixstyle{F}$ the model matrix defined by $\matrixstyle{F}^{(i,j)}=f_i(\vectorstyle{x}_j)$ and $\vectorstyle{\beta}$ the coefficients of the multivariate polynomial trend: $m(\vectorstyle{x})=\vectorstyle{\beta}^\top\vectorstyle{f}(\vectorstyle{x})$.

It is now possible to obtain a predictor of the Gaussian process by assessing $p[\mathcal{y}(\vectorstyle{x})|\matrixstyle{S}^{\vectorstyle{i}}]$, which has a closed expression in the form of a normal distribution of which the mean $\mathbb{E}[\mathcal{y}(\vectorstyle{x}|\matrixstyle{S}^{\vectorstyle{i}})]\triangleq\mu(\vectorstyle{x})$ and variance $\mathbb{V}[\mathcal{y}(\vectorstyle{x}|\matrixstyle{S}^{\vectorstyle{i}})]\triangleq\Sigma(\vectorstyle{x})$ are given by 
\begin{subequations}
	\begin{align}
		\mu(\vectorstyle{x})=&\:\vectorstyle{\beta}^\top\cdot\vectorstyle{f}(\vectorstyle{x})+\vectorstyle{\alpha}^\top\cdot\vectorstyle{\psi}(\vectorstyle{x}) \label{eq:Y}\\ 
		\Sigma(\vectorstyle{x})=&\:\sigma^2\bigl\{\psi(0)-||\vectorstyle{\psi}(\vectorstyle{x})||^2_{\matrixstyle{\Psi}^{-1}}+||\vectorstyle{g}(\vectorstyle{x})||^2_{\matrixstyle{\Gamma}}\bigr\} \label{eq:s2}
	\end{align}
\end{subequations}

\noindent with $\vectorstyle{\psi}(\vectorstyle{x})=[\psi(\vectorstyle{x}_{1},\vectorstyle{x}),...,\psi(\vectorstyle{x}_{n},\vectorstyle{x})]^\top$,
$\vectorstyle{\alpha}=\matrixstyle{\Psi}^{-1}(\vectorstyle{y}-\matrixstyle{F}\vectorstyle{\beta})$, 
$\matrixstyle{\Gamma}=(\matrixstyle{F}^\top\matrixstyle{\Psi}^{-1}\matrixstyle{F})^{-1}$,
$\vectorstyle{g}(\textbf{x})=\matrixstyle{F}^\top\matrixstyle{\Psi}^{-1}\vectorstyle{\psi}(\textbf{x})-\vectorstyle{f}(\vectorstyle{x})$.

In this work, we have made use of the open-source \MATLAB toolbox ooDACE (\textsl{object-oriented Design and Analysis of Computer Experiments}) \cite{couckuyt:2014}. To solve the MLE problem, a genetic algorithm is employed to find the region of the global optimizer, which is subsequently refined by means of a \textsl{sequential quadratic programming} (SQP) approach. 

\subsection{Regressive treatment}

As its name indicates, GPI exactly interpolates between function evaluations, which is a desirable characteristic if these are deterministic. However, when dealing with stochastic calls, for example due to noise, this interpolative character might lead to overfitting and `spiky' behaviour. Therefore, a regressive or stochastic formulation might be preferred. \textsl{Gaussian Process Regression} (GPR) can be realized by adding a `nugget' (inspired by its geostatistical origin) to the covariance matrix. This nugget takes on the form of a noise-to-signal ratio matrix such that the covariance matrix becomes: $\matrixstyle{\Psi}_r\triangleq\matrixstyle{\Psi}+\sigma_r^{-2}\cdot10^{\lambda}\cdot\matrixstyle{I}$, with $\lambda$ a measure of noise. Estimation of this constant can be realized using a MLE as has been introduced before. The polynomial constants and process variance can be updated to respectively give
\begin{subequations}
	\begin{align}
		\vectorstyle{\beta}_r=&\:(\matrixstyle{F}^\top\matrixstyle{\Psi}_r^{-1}\matrixstyle{F})^{-1}\matrixstyle{F}^\top\matrixstyle{\Psi}_r^{-1}\textbf{y} \\ 
		\sigma_r^2=&\:n(\matrixstyle{S}^{\vectorstyle{i}})^{-1}(\vectorstyle{y}-\matrixstyle{F}\vectorstyle{\beta})^\top\matrixstyle{\Psi}_r^{-1}(\vectorstyle{y}-\matrixstyle{F}\vectorstyle{\beta})
	\end{align}
\end{subequations}

The MLE problem now becomes
\begin{equation}\label{eq:MLEr}
	\vectorstyle{\theta}_*,\lambda_*=\arg\max_{\vectorstyle{\theta},\lambda}\: -n(\matrixstyle{S}^{\vectorstyle{i}})\log\sigma_r^2-\log|\matrixstyle{\Psi}_r|
\end{equation}

As was done before, we can build forth on the definition of a Gaussian process and the theorem of Bayes the predictive distribution $p[\mathcal{y}(\vectorstyle{x})|\matrixstyle{S}^{\vectorstyle{i}}]$ can be directly evaluated and gives again a normal distribution of which the mean $\mathbb{E}[\mathcal{y}(\vectorstyle{x}|\matrixstyle{S}^{\vectorstyle{i}})]\triangleq\mu_r(\vectorstyle{x})$ and variance $\mathbb{V}[\mathcal{y}(\vectorstyle{x}|\matrixstyle{S}^{\vectorstyle{i}})]\triangleq\Sigma_r(\vectorstyle{x})$ can be directly evaluated 
\begin{subequations}
	\begin{align}
		\mu_r(\vectorstyle{x})=&\:\vectorstyle{\beta}_r^\top\cdot\vectorstyle{f}(\vectorstyle{x})+\vectorstyle{\alpha}_r^\top\cdot\vectorstyle{\psi}(\vectorstyle{x}) \label{eq:Yr}\\ 
		\Sigma_r(\vectorstyle{x})=&\:\sigma_r^2\bigl\{\psi(0)-||\vectorstyle{\psi}(\vectorstyle{x})||^2_{\matrixstyle{\Psi}_r^{-1}}+||\vectorstyle{g}_r(\vectorstyle{x})||^2_{\matrixstyle{\Gamma}_r}\bigr\} \label{eq:s2r}
	\end{align}
\end{subequations}
\noindent with updated function $\vectorstyle{\alpha}_r=\matrixstyle{\Psi}_r^{-1}(\vectorstyle{y}-\matrixstyle{F}\vectorstyle{\beta}_r)$, 
$\matrixstyle{\Gamma}_r=(\matrixstyle{F}^\top\matrixstyle{\Psi}_r^{-1}\matrixstyle{F})^{-1}$,
$\vectorstyle{g}_r(\textbf{x})=\matrixstyle{F}^\top\matrixstyle{\Psi}_r^{-1}\vectorstyle{\psi}(\textbf{x})-\vectorstyle{f}(\vectorstyle{x})$.

\subsection{Acquisition function}

In the context of optimization, we are specifically interested in evaluating the design that on the one hand is expected to outperform the current best-evaluated point, denoted as $\text{y}_{min}=\min[y(\matrixstyle{X})]$, and on the other hand contribute most strongly to emulation of the objective function by Gaussian process interpolator. This can be realized by maximizing the expected improvement $\mathbb{E}_{ \mathcal{y}}[\mathcal{I}(\vectorstyle{x}|\matrixstyle{S}^{\vectorstyle{i}})]$, with the improvement defined as $\mathcal{I}(\vectorstyle{x}|\text{y}_{min}) = \max[\text{y}_{min}-\mathcal{y}(\vectorstyle{x}),0]$ \cite{mockus:1978}, such that
\begin{align} \label{eq:EI}
	\text{EI}(\vectorstyle{x}|\text{y}_{min})\triangleq&\:\mathbb{E}_{\mathcal{y}}[\mathcal{I}(\vectorstyle{x}|\text{y}_{min})] \nonumber\\ 
	=&\:(y_{min}-\mu(\vectorstyle{x}))\cdot\Phi(y_{min}|\mu(\vectorstyle{x}),\Sigma(\vectorstyle{x})) \nonumber\\ 
	&+\Sigma(\vectorstyle{x})\cdot\phi(y_{min}|\mu(\vectorstyle{x}),\Sigma(\vectorstyle{x})) 
\end{align}

This forms the basis of the \textsl{efficient global optimization} (EGO) algorithm (Algorithm~\ref{alg:EGO}\footnote{$k_{EI}$ is introduced as an additional constant that serves as a stopping criterion: if no further significant improvement can be made, stop ($\texttt{break}$) the optimization.}) \cite{jones:1998}.

\begin{algorithm}[H]
	\caption{Efficient Global Optimization (EGO) \cite{jones:1998}}
	\label{alg:EGO}
	\begin{algorithmic}[1]
		\Require Evaluated sampling plan $\matrixstyle{S}^{\vectorstyle{i}}=\{\matrixstyle{X},\vectorstyle{y}\}$ 
		\While {computational budget is not exhausted}
		\State $\vectorstyle{\theta}_{*}\gets\arg\max_{\vectorstyle{\theta}}\mathcal{L}^{\vectorstyle{i}}(\vectorstyle{\theta}|\matrixstyle{S}^{\vectorstyle{i}})$ \Comment{Equation~\ref{eq:MLE}}
		\State $\text{y}_{min}\gets\min\:\vectorstyle{y}$
		\State $\text{z}_* \gets \max_{\vectorstyle{x}}\:\text{EI}(\vectorstyle{x}|\vectorstyle{\theta}_{*},\text{y}_{min})$
		\If {$\text{z}_*>k_{EI}$}
		\State $\vectorstyle{x}_* = \arg\max_{\vectorstyle{x}}\:\text{EI}(\vectorstyle{x}|\vectorstyle{\theta}_{*},\text{y}_{min})$ \Comment{Equation~\ref{eq:EI}}
		\State $\text{y}_*\gets\text{y}(\vectorstyle{x}_*)$
		\State $\matrixstyle{S}^{\vectorstyle{i}}\gets\matrixstyle{S}^{\vectorstyle{i}}\cup\{\vectorstyle{x}_*,\text{y}_*\}$
		\Else
		\State $\texttt{break}$
		\EndIf
		\EndWhile
		\State $\text{x}_{*}\gets\arg\min_{\vectorstyle{x}}\:\vectorstyle{y}$
	\end{algorithmic}
\end{algorithm}

\subsection{Multi-objective treatment}

When confronted with multiple objectives, we are interested in improving upon the Pareto front. This can be realized by assessing the hypervolume indicator $\mathcal{H}(\mathcal{P})$, which corresponds to the Lebesque measure of the hyperspace dominated by the Pareto front bounded by a reference point in the objective space $\vectorstyle{r}\triangleq\{y_j^{max}+\epsilon|j=1...d(\mathcal{Y})\}$. In a manner similar to the single-objective approach, the indicator can be used to define a scalar improvement function $\mathcal{I}_{hv}(\vectorstyle{x}|\mathcal{P},\vectorstyle{r})$ which measures the contribution (or improvement) of the point $\vectorstyle{y}(\vectorstyle{x})$ to the Pareto set $\mathcal{P}$. To realize this, the exclusive hypervolume indicator $\mathcal{H}^e(\vectorstyle{x}|\mathcal{P},\vectorstyle{r})\triangleq\mathcal{H}(\mathcal{P}\cup\vectorstyle{\mathcal{y}}(\vectorstyle{x})|\vectorstyle{r})-\mathcal{H}(\mathcal{P}|\vectorstyle{r})$ is introduced, such that

\begin{align}
	\mathcal{I}_{hv}(\vectorstyle{x}|\mathcal{P},\vectorstyle{r})=&\:\text{max}[\mathcal{H}^e(\vectorstyle{x}|\mathcal{P},\vectorstyle{r}),0]  
\end{align}
		
The expected hypervolume improvement can now be defined as $\mathbb{E}_{\vectorstyle{\mathcal{y}}}[\mathcal{I}_{hv}(\vectorstyle{x}|\matrixstyle{\Phi}_{\vectorstyle{\mathfrak{x}}},\mathcal{P},\vectorstyle{r})]$. The closed form approximate evaluation of this integral as proposed by Couckuyt et al. has been used in this work \cite{couckuyt:2014b}.
\begin{subequations}
	\begin{align}		
		\text{HEI}(\vectorstyle{x}|\mathcal{P},\vectorstyle{r})\triangleq&\:\mathbb{E}_{\vectorstyle{\mathcal{y}}}[\mathcal{I}_{hv}(\vectorstyle{x}|\mathcal{P},\vectorstyle{r})] \nonumber\\
		=&\:\int\mathcal{I}_{hv}(\vectorstyle{x}|\mathcal{P},\vectorstyle{r})\prod\nolimits_{j=1}^{d(\mathcal{Y})} \phi_j(\vectorstyle{x})\cdot d\mathcal{y}_j
	\end{align}
\end{subequations}

\subsection{Crash Constraint Treatment}\label{sec:crash}

In order to enhance EGO such that function evaluation failures can be accounted for, information on whether ($\text{c}_i=+1$) or not ($\text{c}_i=-1$) the design vector ($\vectorstyle{x}_i$) could be evaluated, is stored. In this work we make use of a discriminative \textsl{Gaussian process classifier} (GPC) to model the feasible space directly $p(\mathcal{c}(\vectorstyle{x})|\matrixstyle{S}^{\vectorstyle{c}})$ with $\matrixstyle{S}^{\vectorstyle{c}}=\{\matrixstyle{X},\vectorstyle{c}\}=\{(\vectorstyle{x}_i,\text{c}_i)|i=1,...,n(\matrixstyle{S}^{\vectorstyle{c}})\}$. The idea behind GPC corresponds to turning the output of a regression model into a class probability using a response function $\lambda(\vectorstyle{z}):\mathbb{R}\mapsto[0,1]$, guaranteeing a valid probabilistic interpretation. 

Assessing the GPC now becomes a two step process, in which first the distribution over the latent output is determined $p(\mathcal{y}(\vectorstyle{x})|\matrixstyle{S}^{\vectorstyle{c}})$, after which this distribution is used to make a prediction of the probability of a class $p(\mathcal{c}(\vectorstyle{x})|\matrixstyle{S}^{\vectorstyle{c}})$. The evaluation of these integrals is intractable and must be approximated. In this work use is made of \textsl{expectation propagation} (EP). The EP algorithm relies on the minimization of the \textsl{Kullback-Leibler} (KL) divergence of the exact posterior over the approximate $\text{KL}(p(\vectorstyle{\mathcal{y}}|\matrixstyle{S}^{\vectorstyle{c}})||q(\vectorstyle{\mathcal{y}}|\matrixstyle{S}^{\vectorstyle{c}}))$. This translates itself in an iterative process during which the approximate distribution is sequentially updated. Details can be found in \cite{rasmussen:2006}, chapter 3.

The availability of this approximate Gaussian posterior distribution obtained by means of the EP algorithm now permits the evaluation of the approximate predictive distribution.

\begin{equation}\label{eq:q}
	q(\mathcal{c}(\vectorstyle{x})=1|\matrixstyle{S}^{\vectorstyle{c}})=\Phi\left(\frac{\vectorstyle{\psi}(\vectorstyle{x})^\top\hat{\matrixstyle{\Sigma}}^{-1}\tilde{\vectorstyle{\mu}}}{\sqrt{1+\psi(\vectorstyle{x},\vectorstyle{x})-||\vectorstyle{\psi}(\vectorstyle{x})||^2_{\hat{\matrixstyle{\Sigma}}^{-1}}}}\right)
\end{equation}
\noindent with $\tilde{\matrixstyle{\Sigma}}=\texttt{diag}(\tilde{\sigma}_i^2)$, $\vectorstyle{\mu}=\matrixstyle{\Sigma}^{-1}\tilde{\matrixstyle{\Sigma}}\tilde{\vectorstyle{\mu}}$, $\matrixstyle{\Sigma}=(\matrixstyle{\Psi}^{-1}+\tilde{\matrixstyle{\Sigma}}^{-1})$, $\hat{\matrixstyle{\Sigma}}=\matrixstyle{\Psi}+\tilde{\matrixstyle{\Sigma}}$ and $k_{EP}=q(\matrixstyle{S}^{\vectorstyle{c}})$ determined by means of the EP algorithm.

To obtain the hyperparameters that fully describe the GPC, a MLE problem is solved, as was done for the GPI. The marginal log likelihood can now be approximated as $\log\:\mathcal{L}^{\vectorstyle{c}}(\vectorstyle{\theta}|\matrixstyle{S}^{\vectorstyle{c}})\triangleq\log\:p(\matrixstyle{S}^{\vectorstyle{c}}|\vectorstyle{\theta})\approx\log\:q(\matrixstyle{S}^{\vectorstyle{c}}|\vectorstyle{\theta})$ such that
\begin{align}\label{eq:MLEclassifier}
	\vectorstyle{\theta}_*=&\arg\max_{\vectorstyle{\theta}}\:-\log|\hat{\matrixstyle{\Sigma}}|-\tilde{\vectorstyle{\mu}}^\top\hat{\matrixstyle{\Sigma}}^{-1}\tilde{\vectorstyle{\mu}} \nonumber\\
	&+\sum_{i=1}^{n(\matrixstyle{S}^{\vectorstyle{c}})}\left\{2\log\left[\Phi\left(\frac{\text{c}_i\mu_{\tilde i}}{\sqrt{1+\sigma_{\sim i}^2}}\right)\hat{\sigma}_i\right]+\frac{\hat{\mu}_i^2}{\hat{\sigma}^2_i}\right\}
\end{align}
\noindent with $\hat{\sigma}^2_i\triangleq\sigma^2_{\sim i}+\tilde{\sigma}^2_i$ and $\hat{\mu}_i\triangleq\mu_{\sim i}-\tilde{\mu}_i$.

The GPC can now be used to enhance the acquisition function in such a way that only the feasible domain is searched during its optimization. This is realized by taking the product of the acquisition function and the predictive classifier. The result is a \textsl{`safe' expected improvement} (s-EI):
\begin{equation}
	\begin{aligned}
		\text{s-EI}(\vectorstyle{x}|\vectorstyle{y}_{min},\matrixstyle{S}^{\vectorstyle{c}})\triangleq&\:\mathbb{E}_{\mathcal{y}}[\mathcal{I}(\vectorstyle{x}|\text{y}_{min})]\cdot q(\mathcal{c}(\vectorstyle{x})=1|\matrixstyle{S}^{\vectorstyle{c}})
	\end{aligned}
\end{equation}

The ooDACE toolbox is extended with the publicly available implementation of Rasmussen \& Williams. The MLE problem is solved using a conjugate gradient approach \cite{rasmussen:2006}. 

\section*{Nomenclature}
\subsection*{Acronyms}
\addcontentsline{toc}{section}{Nomenclature}
\begin{IEEEdescription}[\IEEEusemathlabelsep\IEEEsetlabelwidth{$V_1,V_2,V_3$}]	
	\item[AVL]{Athena Vortex Lattice}
	\item[AR]{Aspect Ratio}
	\item[BO]{Bayesian Optimization}
	\item[BWB]{Blended Wing Body}
	\item[CDP]{Conceptual Design Problem}
	\item[CLD]{Closed Loop Dynamics}
	\item[DoE]{Design of Experiments}
	\item[FCP]{Feedback Control Problem}
	\item[(s-)(H)EI]{(Safe-)(Hypervolume) Expected Improvement}
	\item[GPI/C]{Gaussian Process Interpolator/Classifier}
	\item[L/HF(M)]{Low/High Fidelity (Model)}
	\item[MDO]{Multidisciplinary Design Optimization}
	\item[MLE]{Maximum Likelihood Estimation}
	\item[NLP]{Non-Linear Problem}
	\item[OCP]{Optimal Control Problem}
	\item[SQM]{Sequential Quadratic Programming}
	\item[UAS]{Unmanned Aerial System}
	\item[UP]{Uncertainty Propagation}
	\item[XDSM]{Extended Design System Matrix}
\end{IEEEdescription}
	
\subsection*{Greek symbols}
\begin{IEEEdescription}[\IEEEusemathlabelsep\IEEEsetlabelwidth{$V_1,V_2,V_3$}]	
	\item[$\alpha$]{Angle of attack}
	\item[$\vectorstyle{\alpha}$]{GP predictive mean coefficients}
	\item[$\beta$]{Sideslip angle}
	\item[$\vectorstyle{\beta}$]{Polynomial coefficients}
	\item[$\vectorstyle{\epsilon}$]{Irreducible stochasticity representation}
	\item[$\delta$]{Control surface deflection}
	\item[$\gamma$]{Wing twist}
	\item[$\vectorstyle{\gamma}$]{State offset integration}
	\item[$\lambda$]{Regressive (nugget) GP parameter}
	\item[$\vectorstyle{\Omega}$]{Rotational velocity vector}
	\item[$\vectorstyle{\upsilon}$]{Control signal vector}
	\item[$\matrixstyle{\Upsilon}$]{Control signal space}
	\item[$\phi$]{Roll angle}
	\item[$\psi$]{Yaw angle}
	\item[$\sigma^2$]{Process variance}
	\item[$\vectorstyle{\tau}$]{Moment vector}
	\item[$\theta$]{Pitch angle}
	\item[$\vectorstyle{\theta}$]{Covariance function hyperparameters}
	\item[$\vectorstyle{\xi}$]{System state vector}	
	\item[$\vectorstyle{\Xi}$]{System state space}	
\end{IEEEdescription}
	
\subsection*{Roman symbols}
\begin{IEEEdescription}[\IEEEusemathlabelsep\IEEEsetlabelwidth{$V_1,V_2,V_3$}]	
	\item[$b$]{Span width}
	\item[$\vectorstyle{c}$]{B-spline weighting parameters}
	\item[$c$]{Chord length}
	\item[$C_{\{\cdot\}}$]{$\pi$-motivated force/moment coefficient}
	\item[$d(\cdot)$]{Dimensionality}
	\item[$\vectorstyle{d}$]{Conceptual design vector}
	\item[$\matrixstyle{D}$]{Conceptual design set}
	\item[$\mathcal{D}$]{Conceptual design space}
	\item[$e$]{Span efficiency factor}
	\item[D]{Drag force}
	\item[$f_{\text{con}}$]{Control fraction}
	\item[$\vectorstyle{f}$]{Force vector}
	\item[$f_a$]{Acquisition function value}
	\item[$\vectorstyle{g}$]{Gravitational force vector}
	\item[$\mathcal{H}$]{Hypervolume indicator}
	\item[$\matrixstyle{I}$]{Inertia matrix}
	\item[$\mathcal{I}$]{Improvement indicator}
	\item[$k$]{Indiced drag coefficient}
	\item[$\vectorstyle{k}$]{Control parameter vector}
	\item[$\matrixstyle{K}$]{Control parameter set}
	\item[$K_{\{\cdot\}}$]{$\phi$-motivated force/moment coefficient}
	\item[l]{Time integrand objective function}
	\item[L]{Lift force}
	\item[$m$]{Mass}
	\item[$n(\cdot)$]{Set size}
	\item[$\vectorstyle{p}$]{Position vector}
	\item[$\vectorstyle{q}$]{Attitude vector}
	\item[$\mathcal{P}$]{Pareto front}
	\item[$\matrixstyle{S}$]{Training set}
	\item[$S$]{Wetted area}
	\item[$t$]{Time instance}
	\item[$T$]{Mission duration}
	\item[$\mathcal{T}$]{Time grid}
	\item[T]{Thrust}
	\item[$\vectorstyle{v}$]{Velocity vector}
	\item[$w$]{Weight coefficient}
\end{IEEEdescription}
	
\subsection*{Operators}
\begin{IEEEdescription}[\IEEEusemathlabelsep\IEEEsetlabelwidth{$V_1,V_2,V_3$}]	
	\item[$\mathbb{E}_{\{\cdot\}}{[\cdot]}$]{Expected value operator}
	\item[$\vectorstyle{f}(\cdot)$]{Trend vector function}
	\item[$\mathcal{F}_{\{\cdot\}}{[\cdot]}$]{State space dynamics}
	\item[$\vectorstyle{g}(\cdot)$]{Inequality constraint function}	
	\item[$\vectorstyle{h}(\cdot)$]{Equality constraint function}
	\item[$\mathcal{H}{[\cdot]}$]{Hypervolume measure}
	\item[$\mathcal{L}(\cdot)$]{Likelihood function}	
	\item[$\mu_{\{\cdot\}}(\cdot)$]{Predictive mean}
	\item[$\mathcal{N}(\cdot,\cdot)$]{Normal distribution}
	\item[$p(\cdot)$]{Probability density function}
	\item[$q(\cdot)$]{Approximate probability density function}
	\item[$\vectorstyle{\pi}(\cdot)$]{Control law}
	\item[$\phi(\cdot)$]{Standard normal probability density function}
	\item[$\Phi(\cdot)$]{Standard normal cumulative distribution function}
	\item[$\psi(\cdot,\cdot)$]{Covariance function}
	\item[$\Sigma_{\{\cdot\}}(\cdot)$]{Predictive variance}
	\item[$\mathbb{V}_{\{\cdot\}}{[\cdot]}$]{Variance operator}
	\item[$\mu_{\{\cdot\}}(\cdot)$]{Predictive mean}
	\item[$\Upsilon{[\cdot]}$]{Differential flatness transform for $\vectorstyle{\upsilon}$}
	\item[$\Xi{[\cdot]}$]{Differential flatness transform for $\vectorstyle{\xi}$}
	\item[$\vectorstyle{y}(\cdot)$]{Objective vector function}
	\item[$\mathcal{y}(\cdot)$]{Gaussian process}
	
\end{IEEEdescription}
	
\subsection*{Subscript}
\begin{IEEEdescription}[\IEEEusemathlabelsep\IEEEsetlabelwidth{$V_1,V_2,V_3$}]	
	\item[$*$]{Optimal/Reference value}
	\item[$c$]{Control}
	\item[cog]{Center of gravity}
	\item[$hv$]{Hypervolume}
	\item[$r$]{Regressive}
\end{IEEEdescription}
	
\subsection*{Superscript}
\begin{IEEEdescription}[\IEEEusemathlabelsep\IEEEsetlabelwidth{$V_1,V_2,V_3$}]
	\item[$\vectorstyle{c}$]{Classifier}
	\item[$\vectorstyle{i}$]{Interpolator}
	\item[global]{Global inertial frame of reference}
	\item[local]{Body fixed frame of reference}
	\item[$\top$]{Transpose}
\end{IEEEdescription}

\bibliographystyle{IEEEtranN}
{\small\bibliography{Wauters.bib}}

\begin{thebibliography}{46}
\providecommand{\natexlab}[1]{#1}
\providecommand{\url}[1]{#1}
\csname url@samestyle\endcsname
\providecommand{\newblock}{\relax}
\providecommand{\bibinfo}[2]{#2}
\providecommand{\BIBentrySTDinterwordspacing}{\spaceskip=0pt\relax}
\providecommand{\BIBentryALTinterwordstretchfactor}{4}
\providecommand{\BIBentryALTinterwordspacing}{\spaceskip=\fontdimen2\font plus
\BIBentryALTinterwordstretchfactor\fontdimen3\font minus
  \fontdimen4\font\relax}
\providecommand{\BIBforeignlanguage}[2]{{%
\expandafter\ifx\csname l@#1\endcsname\relax
\typeout{** WARNING: IEEEtranN.bst: No hyphenation pattern has been}%
\typeout{** loaded for the language `#1'. Using the pattern for}%
\typeout{** the default language instead.}%
\else
\language=\csname l@#1\endcsname
\fi
#2}}
\providecommand{\BIBdecl}{\relax}
\BIBdecl

\bibitem[Saeed et~al.(2018)Saeed, Younes, Cai, and Cai]{saeed:2018}
A.~S. Saeed, A.~B. Younes, C.~Cai, and G.~Cai, ``A survey of hybrid unmanned
  aerial vehicles,'' \emph{Progress in Aerospace Sciences}, vol.~98, pp.
  91--105, 2018.

\bibitem[Ducard and Allenspach(2021)]{ducard:2021}
G.~J.~J. Ducard and M.~Allenspach, ``Review of designs and flight control
  techniques of hybrid and convertible {VTOL} {UAVs},'' \emph{Aerospace Science
  and Technology}, vol. 118, p. 107035, 2021.

\bibitem[Smeur et~al.(2019)Smeur, Bronz, and de~Croon]{smeur:2019}
E.~J.~J. Smeur, M.~Bronz, and G.~C. H.~E. de~Croon, ``Incremental control and
  guidance of hybrid aircraft applied to a tailsitter unmanned air vehicle,''
  \emph{Journal of Guidance, Control, and Dynamics}, vol.~43, no.~2, pp.
  274--287, 2019.

\bibitem[Fuhrer et~al.(2019)Fuhrer, Verling, Stastny, and
  Siegwart]{fuhrer:2019}
S.~Fuhrer, S.~Verling, T.~Stastny, and R.~Siegwart, ``Fault-tolerant flight
  control of a {VTOL} tailsitter {UAV},'' in \emph{2019 International
  Conference on Robotics and Automation (ICRA)}, 2019, Conference Proceedings,
  pp. 4134--4140.

\bibitem[Olsson et~al.(2021)Olsson, Verling, Stastny, and
  Siegwart]{olsson:2021}
C.~Olsson, S.~L. Verling, T.~Stastny, and R.~Siegwart, \emph{Full Envelope
  System Identification of a {VTOL} Tailsitter {UAV}}, ser. AIAA SciTech
  Forum.\hskip 1em plus 0.5em minus 0.4em\relax American Institute of
  Aeronautics and Astronautics, 2021.

\bibitem[Theys and De~Schutter(2016)]{theys:2016}
B.~Theys and J.~De~Schutter, ``Parameter selection method and performance
  assessment for the preliminary design of electrically powered transitioning
  {VTOL} {UAVs},'' in \emph{Proceedings of IMAV 2016}, 2016, Conference
  Proceedings, pp. 221--228.

\bibitem[Banazadeh and Taymourtash(2015)]{banazadeh:2015}
\BIBentryALTinterwordspacing
A.~Banazadeh and N.~Taymourtash, ``Optimal control of an aerial tail sitter in
  transition flight phases,'' \emph{Journal of Aircraft}, vol.~53, no.~4, pp.
  914--921, 2015. [Online]. Available: \url{https://doi.org/10.2514/1.C033339}
\BIBentrySTDinterwordspacing

\bibitem[Oosedo et~al.(2017)Oosedo, Abiko, Konno, and Uchiyama]{oosedo:2017}
A.~Oosedo, S.~Abiko, A.~Konno, and M.~Uchiyama, ``Optimal transition from
  hovering to level-flight of a quadrotor tail-sitter {UAV},'' \emph{Autonomous
  Robots}, vol.~41, no.~5, pp. 1143--1159, 2017.

\bibitem[Verling et~al.(2017)Verling, Stastny, Bättig, Alexis, and
  Siegwart]{verling:2017}
S.~Verling, T.~Stastny, G.~Bättig, K.~Alexis, and R.~Siegwart, ``Model-based
  transition optimization for a {VTOL} tailsitter,'' in \emph{2017 IEEE
  International Conference on Robotics and Automation (ICRA)}, 2017, pp.
  3939--3944.

\bibitem[Li et~al.(2020)Li, Sun, Zhou, Wen, Low, and Chen]{li:2020e}
B.~Li, J.~Sun, W.~Zhou, C.~Y. Wen, K.~H. Low, and C.~K. Chen, ``Transition
  optimization for a {VTOL} tail-sitter {UAV},'' \emph{IEEE/ASME Transactions
  on Mechatronics}, vol.~25, no.~5, pp. 2534--2545, 2020.

\bibitem[Garcia-Sanz(2019)]{garcia:2019}
M.~Garcia-Sanz, ``Control co-design: An engineering game changer,''
  \emph{Advanced Control for Applications}, vol.~1, no.~1, p. e18, 2019.

\bibitem[Herber and Allison(2018)]{herber:2018}
D.~R. Herber and J.~T. Allison, ``Nested and simultaneous solution strategies
  for general combined plant and control design problems,'' \emph{Journal of
  Mechanical Design}, vol. 141, no.~1, 2018.

\bibitem[Kaneko and Martins(2023)]{kaneko:2023}
S.~Kaneko and J.~R. Martins, ``Simultaneous optimization of conceptual design
  and takeoff trajectory of a lift-plus-cruise {UAV},'' in \emph{10th
  Autonomous {VTOL} Technical Meeting}, 2023, Conference Proceedings.

\bibitem[Wauters et~al.(2022)Wauters, Lefebvre, and Crevecoeur]{wauters:2022c}
J.~Wauters, T.~Lefebvre, and G.~Crevecoeur, ``Comparative study of co-design
  strategies for mission-specific design of quadcopters using differential
  flatness and bayesian optimization,'' in \emph{2022 IEEE/ASME International
  Conference on Advanced Intelligent Mechatronics (AIM)}.\hskip 1em plus 0.5em
  minus 0.4em\relax IEEE, 2022, Conference Proceedings, pp. 703--709.

\bibitem[Wauters et~al.(2023)Wauters, Lefebvre, and Crevecoeur]{wauters:2023}
------, ``Multi-objective co-design for mission-specific development of
  unmanned aerial systems,'' in \emph{2023 IEEE/ASME International Conference
  on Advanced Intelligent Mechatronics (AIM)}, 2023, Conference Proceedings,
  pp. 17--24.

\bibitem[Balesdent et~al.(2022)Balesdent, Brevault, Valderrama-Zapata, and
  Urbano]{balesdent:2022}
M.~Balesdent, L.~Brevault, J.-L. Valderrama-Zapata, and A.~Urbano,
  ``All-at-once formulation integrating pseudo-spectral optimal control for
  launch vehicle design and uncertainty quantification,'' \emph{Acta
  Astronautica}, vol. 200, pp. 462--477, 2022.

\bibitem[Morita et~al.(2020)Morita, Tsuchiya, and Taguchi]{morita:2020}
N.~Morita, T.~Tsuchiya, and H.~Taguchi, \emph{MDO of Hypersonic Waverider with
  Trajectory-Aero-Structure Coupling}, ser. International Space Planes and
  Hypersonic Systems and Technologies Conferences.\hskip 1em plus 0.5em minus
  0.4em\relax American Institute of Aeronautics and Astronautics, 2020.

\bibitem[Matos and Marta(2022)]{matos:2022}
N.~M.~B. Matos and A.~C. Marta, ``Concurrent trajectory optimization and
  aircraft design for the air cargo challenge competition,'' \emph{Aerospace},
  vol.~9, no.~7, 2022.

\bibitem[Lupp et~al.(2022)Lupp, Clark, Aksland, and Alleyne]{lupp:2022}
C.~A. Lupp, D.~L. Clark, C.~T. Aksland, and A.~G. Alleyne, \emph{Mission and
  Shape Optimization of a HALE Aircraft including Transient Power and Thermal
  Constraints}, ser. AIAA AVIATION Forum.\hskip 1em plus 0.5em minus
  0.4em\relax American Institute of Aeronautics and Astronautics, 2022.

\bibitem[Hendricks et~al.(2020)Hendricks, Aretskin-Hariton, Ingraham, Gray,
  Schnulo, Chin, Falck, and Hall]{hendricks:2020}
E.~S. Hendricks, E.~Aretskin-Hariton, D.~Ingraham, J.~S. Gray, S.~L. Schnulo,
  J.~Chin, R.~Falck, and D.~Hall, \emph{Multidisciplinary Optimization of an
  Electric Quadrotor Urban Air Mobility Aircraft}, ser. AIAA AVIATION
  Forum.\hskip 1em plus 0.5em minus 0.4em\relax American Institute of
  Aeronautics and Astronautics, 2020.

\bibitem[Delbecq et~al.(2021)Delbecq, Budinger, Coic, and
  Bartoli]{delbecq:2021}
S.~Delbecq, M.~Budinger, C.~Coic, and N.~Bartoli, \emph{Trajectory and design
  optimization of multirotor drones with system simulation}, ser. AIAA SciTech
  Forum.\hskip 1em plus 0.5em minus 0.4em\relax American Institute of
  Aeronautics and Astronautics, 2021.

\bibitem[Martins and Lambe(2013)]{martins:2013}
J.~R. R.~A. Martins and A.~B. Lambe, ``Multidisciplinary design optimization: A
  survey of architectures,'' \emph{AIAA Journal}, vol.~51, no.~9, pp.
  2049--2075, 2013.

\bibitem[Ha et~al.(2018)Ha, Coros, Alspach, Kim, and Yamane]{ha:2018}
S.~Ha, S.~Coros, A.~Alspach, J.~Kim, and K.~Yamane, ``Computational
  co-optimization of design parameters and motion trajectories for robotic
  systems,'' \emph{The International Journal of Robotics Research}, vol.~37,
  no. 13-14, pp. 1521--1536, 2018.

\bibitem[Whitman and Choset(2018)]{whitman:2018}
J.~Whitman and H.~Choset, ``Task-specific manipulator design and trajectory
  synthesis,'' \emph{IEEE Robotics and Automation Letters}, vol.~4, no.~2, pp.
  301--308, 2018.

\bibitem[Toussaint et~al.(2021)Toussaint, Ha, and Oguz]{toussaint:2021}
M.~Toussaint, J.-S. Ha, and O.~S. Oguz, ``Co-optimizing robot, environment, and
  tool design via joint manipulation planning,'' in \emph{2021 IEEE
  International Conference on Robotics and Automation (ICRA)}.\hskip 1em plus
  0.5em minus 0.4em\relax IEEE, 2021, pp. 6600--6606.

\bibitem[Lefebvre et~al.(2023)Lefebvre, Wauters, Ostyn, and
  Crevecoeur]{lefebvre:2023}
T.~Lefebvre, J.~Wauters, F.~Ostyn, and G.~Crevecoeur, ``Towards task tailored
  articulated robot designs,'' in \emph{IEEE/ASME (AIM) International
  Conference on Advanced Intelligent Mechatronics}.\hskip 1em plus 0.5em minus
  0.4em\relax IEEE, 2023, Conference Proceedings.

\bibitem[Tal et~al.(2023)Tal, Ryou, and Karaman]{tal:2023}
E.~Tal, G.~Ryou, and S.~Karaman, ``Aerobatic trajectory generation for a {VTOL}
  fixed-wing aircraft using differential flatness,'' \emph{IEEE Transactions on
  Robotics}, pp. 1--15, 2023.

\bibitem[Ritz and D'Andrea(2017)]{ritz:2017}
R.~Ritz and R.~D'Andrea, ``A global controller for flying wing tailsitter
  vehicles,'' in \emph{2017 IEEE international conference on robotics and
  automation (ICRA)}.\hskip 1em plus 0.5em minus 0.4em\relax IEEE, 2017, pp.
  2731--2738.

\bibitem[Wauters(2024)]{wauters:2024}
\BIBentryALTinterwordspacing
J.~Wauters, ``Ergo-ii: An improved bayesian optimization technique for robust
  design with multiple objectives, failed evaluations, and stochastic
  parameters,'' \emph{Journal of Mechanical Design}, vol. 146, no.~10, 2024.
  [Online]. Available: \url{https://doi.org/10.1115/1.4064674}
\BIBentrySTDinterwordspacing

\bibitem[Drela and Youngren(2017)]{drela:2017}
\BIBentryALTinterwordspacing
M.~Drela and H.~Youngren, ``{AVL} 3.27,'' 2017. [Online]. Available:
  \url{http://web.mit.edu/drela/Public/web/avl/}
\BIBentrySTDinterwordspacing

\bibitem[Drela(2014)]{drela:2014}
M.~Drela, \emph{Flight Vehicle Aerodynamics}.\hskip 1em plus 0.5em minus
  0.4em\relax Cambridge, Massachusetts: The MIT Press, 2014.

\bibitem[Lustosa et~al.(2018)Lustosa, Defaÿ, and Moschetta]{lustosa:2018}
L.~R. Lustosa, F.~Defaÿ, and J.-M. Moschetta, ``Global singularity-free
  aerodynamic model for algorithmic flight control of tail sitters,''
  \emph{Journal of Guidance, Control, and Dynamics}, vol.~42, no.~2, pp.
  303--316, 2018.

\bibitem[Stoical et~al.(2016)Stoical, Iv{\u{a}}nu{\c{s}}c{\u{a}}, Prodan, and
  Popescu]{stoical:2016}
F.~Stoical, V.-M. Iv{\u{a}}nu{\c{s}}c{\u{a}}, I.~Prodan, and D.~Popescu,
  ``Obstacle avoidance via b-spline parametrizations of flat trajectories,'' in
  \emph{2016 24th Mediterranean Conference on Control and Automation
  (MED)}.\hskip 1em plus 0.5em minus 0.4em\relax IEEE, 2016, pp. 1002--1007.

\bibitem[Dubins(1957)]{dubins:1957}
L.~E. Dubins, ``On curves of minimal length with a constraint on average
  curvature, and with prescribed initial and terminal positions and tangents,''
  \emph{American Journal of Mathematics}, vol.~79, no.~3, pp. 497--516, 1957.

\bibitem[Owen et~al.(2014)Owen, Beard, and McLain]{owen:2014b}
M.~Owen, R.~W. Beard, and T.~McLain, \emph{Implementing Dubins airplane paths
  on fixed-wing {UAVs}}.\hskip 1em plus 0.5em minus 0.4em\relax Springer, 2014,
  book section~68.

\bibitem[Williamson et~al.(2012)Williamson, McGranahan, Broughton, Deters,
  Brandt, and Selig]{williamson:2012}
G.~A. Williamson, B.~D. McGranahan, B.~A. Broughton, R.~W. Deters, J.~B.
  Brandt, and M.~S. Selig, \emph{Summary of Low-Speed Airfoil Data}.\hskip 1em
  plus 0.5em minus 0.4em\relax University of Illinois at Urbana-Champaign, IL,
  USA: Department of Aerospace Engineering, 2012, vol.~5.

\bibitem[S\'{o}bester and Forrester(2014)]{sobester:2014}
A.~S\'{o}bester and A.~Forrester, \emph{Aircraft Aerodynamic Design: Geometry
  and Optimization}.\hskip 1em plus 0.5em minus 0.4em\relax Wiley, 2014.

\bibitem[Wauters et~al.(2020)Wauters, Couckuyt, Knudde, Dhaene, and
  Degroote]{wauters:2020b}
J.~Wauters, I.~Couckuyt, N.~Knudde, T.~Dhaene, and J.~Degroote,
  ``Multi-objective optimization of a wing fence on a {UAV} using
  surrogate-derived gradients,'' \emph{Structural and Multidisciplinary
  Optimization}, vol.~61, pp. 353--364, 2020.

\bibitem[Rasmussen and Williams(2006)]{rasmussen:2006}
C.~E. Rasmussen and C.~K.~I. Williams, \emph{Gaussian Processes for Machine
  Learning}.\hskip 1em plus 0.5em minus 0.4em\relax the MIT Press, 2006.

\bibitem[Shahriari et~al.(2016)Shahriari, Swersky, Wang, Adams, and
  Freitas]{shahriari:2016}
B.~Shahriari, K.~Swersky, Z.~Wang, R.~P. Adams, and N.~d. Freitas, ``Taking the
  human out of the loop: A review of bayesian optimization,'' \emph{Proceedings
  of the IEEE}, vol. 104, no.~1, pp. 148--175, 2016.

\bibitem[Keane(2006)]{keane:2006}
A.~J. Keane, ``Statistical improvement criteria for use in multiobjective
  design optimization,'' \emph{AIAA Journal}, vol.~44, no.~4, pp. 879--891,
  2006.

\bibitem[Palar et~al.(2020)Palar, Zuhal, Chugh, and Rahat]{palar:2020}
P.~S. Palar, L.~R. Zuhal, T.~Chugh, and A.~Rahat, \emph{On the Impact of
  Covariance Functions in Multi-Objective Bayesian Optimization for Engineering
  Design}, ser. AIAA SciTech Forum.\hskip 1em plus 0.5em minus 0.4em\relax
  American Institute of Aeronautics and Astronautics, 2020, pp. 1--13.

\bibitem[Couckuyt et~al.(2014{\natexlab{a}})Couckuyt, Dhaene, and
  Demeester]{couckuyt:2014}
I.~Couckuyt, T.~Dhaene, and P.~Demeester, ``oo{DACE} toolbox: a flexible
  object-oriented kriging implementation,'' \emph{J. Mach. Learn. Res.},
  vol.~15, no.~1, pp. 3183--3186, 2014.

\bibitem[Mockus et~al.(1978)Mockus, Tiesis, and Zilinskas]{mockus:1978}
J.~Mockus, V.~Tiesis, and A.~Zilinskas, ``The application of bayesian methods
  for seeking the extremum,'' in \emph{Towards Global Optimization 2:
  Proceedings of a Workshop at the University of Cagliari, Italy, October
  1974}, L.~D. 2 and G.~Szego, Eds., vol.~2, 1978, Conference Proceedings, pp.
  117--129.

\bibitem[Jones et~al.(1998)Jones, Schonlau, and Welch]{jones:1998}
D.~R. Jones, M.~Schonlau, and W.~J. Welch, ``Efficient global optimization of
  expensive black-box functions,'' \emph{Journal of Global Optimization},
  vol.~13, no.~4, pp. 455--492, 1998.

\bibitem[Couckuyt et~al.(2014{\natexlab{b}})Couckuyt, Deschrijver, and
  Dhaene]{couckuyt:2014b}
I.~Couckuyt, D.~Deschrijver, and T.~Dhaene, ``Fast calculation of
  multiobjective probability of improvement and expected improvement criteria
  for pareto optimization,'' \emph{Journal of Global Optimization}, vol.~60,
  no.~3, pp. 575--594, 2014.

\end{thebibliography}

\end{document}